\documentclass[11pt, a4paper,english]{article}
\usepackage{csquotes,amsmath,amsfonts,comment,caption,subcaption,graphicx,url,color,arydshln,mathtools,authblk}
\usepackage[doi=false,isbn=false,url=false,eprint=false,maxnames=50,date=year,clearlang,giveninits=true]{biblatex}
\usepackage{hyperref}
\renewbibmacro{in:}{,}
\usepackage{amsthm}
\usepackage[plain]{fancyref}
\usepackage[a4paper]{geometry}
\usepackage{overpic}
\graphicspath{{./figures/}}
\DeclareFieldFormat[article, book, inbook, incollection, inproceedings, misc, thesis, unpublished]{title}{\textit{#1}}    

\DeclareFieldFormat[article]{journaltitle}{\textnormal{#1}}
\DeclareFieldFormat[incollection]{booktitle}{\textnormal{#1}}

\theoremstyle{remark}
\newtheorem*{remark}{Remark}

\definecolor{ucyan}{RGB}{4, 255, 255}
\newcommand{\ucyan}{\textcolor{ucyan}}

\definecolor{ublue}{RGB}{23, 0, 255}
\newcommand{\ublue}{\textcolor{ublue}}

\definecolor{upurple}{RGB}{128, 0, 128}
\newcommand{\upurple}{\textcolor{upurple}}

\definecolor{ured}{RGB}{255, 0, 0}
\newcommand{\ured}{\textcolor{ured}}

\definecolor{ugreen}{RGB}{0, 255, 0}

\definecolor{uorange}{RGB}{254, 128, 1}

\title{Spline Shallow Water Moment Equations }
\author[1]{Ullika Scholz\footnote{Corresponding author, email address {scholz@acom.rwth-aachen.de}}}
\author[2,3]{Julian Koellermeier}
\affil[1]{Applied and Computational Mathematics, RWTH Aachen University, Aachen, Germany}
\affil[2]{Department of Mathematics, Computer Science and Statistics, Ghent University, Ghent, Belgium}
\affil[3]{Bernoulli Institute, University of Groningen, Groningen, the Netherlands}

\date{\today}

\begin{document}
\maketitle
\begin{abstract}
\noindent
    Reduced models for free-surface flows are required due to the high dimensionality of the underlying incompressible Navier-Stokes equations, which need to fully resolve the flow in vertical direction to compute the surface height.
    On the other hand, standard reduced models, such as the classical Shallow Water Equations (SWE), which assume a small depth-to-length ratio and use depth-averaging, do not provide information about the vertical velocity profile variations. As a compromise, a recently proposed moment approach for shallow flow using Legendre polynomials as ansatz functions for vertical velocity variations showed the derivation of so-called Shallow Water Moment Equations (SWME) that combine low dimensionality with velocity profile modeling. However, only global polynomials are considered so far.
    
    This paper introduces Spline Shallow Water Moment Equations (SSWME) where piecewise defined spline ansatz functions allow for a flexible representation of velocity profiles with lower regularity. The local support of the spline basis functions opens up the possibility of adaptability and greater flexibility regarding some typical profile shapes. We systematically derive and analyze hierarchies of SSWME models with different number of basis functions and different degrees, before deriving a regularized hyperbolic version by performing a hyperbolic regularization with analytical proof of hyperbolicity for a hierarchy of high-order SSWME models. Numerical simulations show high accuracy and robustness of the new models.
\end{abstract}


\section{Introduction}
\label{sec:Introduction}
The Shallow Water Equations (SWE) are a well-studied
model for shallow free surface flows where the fluid height is small in relation to a typical wavelength.
Applications are the prediction 
of tsunamis, floodings, and other geophysical hazards as listed in \cite{KowalskiTorrilhon2019}.
The equations are typically derived from the incompressible Euler equations via vertical averaging, which corresponds to assuming a constant vertical profile in the horizontal velocity \cite{Godlewski2021}. 
This vertical averaging reduces the dimension as the equations no longer depend on the vertical variable.

For certain applications, however, the assumption of a constant velocity profile is too restrictive, even when the shallowness assumption is fulfilled. 
Examples are granular flows or mud flows, where friction forces slow down the fluid at the bottom, which results in large differences in the horizontal flow velocity along the vertical axis \cite{christen2010ramms,kern2009measured,nagl2020velocity,schaefer2013velocity}. Consequentially, the SWE become inaccurate.

To efficiently model such flows without fully resolving the vertical axis, the recently proposed moment method \cite{KowalskiTorrilhon2019} uses a global polynomial expansion of the velocity profile in vertical direction to derive the \emph{Shallow Water Moment Equations} (SWME) hierarchy. The SWME are a coupled PDE system for the coefficients of the polynomial expansion. 
The SWE are included in this hierarchy as the simplest zero-order model. Higher-order SWME models incorporating a larger number of equations have been shown to capture complex profile shapes \cite{KowalskiTorrilhon2019}.

The idea has gained significant interest in recent years:
Non-hydrostatic versions of the SWME can capture dispersive effects \cite{escalante2024,scholz2023} like the  non-hydrostatic layer-wise approaches \cite{Escalante2019,FernandezNieto2018}; hyperbolic corrections, the so-called \emph{Hyperbolic Shallow Water Moment Equations} (HSWME) \cite{KoellermeierRominger2020} guarantee global hyperbolicity; a combination with a layer-wise approach \cite{garresdiaz2023} allows for even more complex shapes. 
Different models were further analyzed and applied \cite{garres2021,Koellermeier2020Eq,Koellermeier2020steady,Steldermann2024}. Extensions to two-dimensional flows have recently appeared \cite{bauerle2024rotationalinvariancehyperbolicityshallow, VerbiestKoellermeier}. 
While several moment models have been derived, they are exclusively limited to a global polynomial Legendre basis expansion of the velocity profile. This limits the potential expressivity and accuracy of the models, as it requires high regularity. It is clear that other basis functions should be explored. However, it remains to be investigated which basis functions are good candidates to start with.

While the moment model approach is relatively new in the field of free-surface flows, moment methods have a long history in the field of
rarefied gas dynamics. This started by the work of Grad \cite{Grad1949}, who uses a weighted Hermite polynomial expansion due to the orthogonality with respect to the Maxwellian equilibrium distribution - similar to the Legendre polynomial basis for free-surface flows.
For gas dynamics, the usefulness of piecewise defined spline basis functions for moment models was demonstrated in \cite{koellermeier2020spline}. It showed that accurate models with beneficial analytical properties like conservation and hyperbolicity can be derived using piecewise defined splines with local support.

This paper derives the first Spline Shallow Water Moment Equations (SSWME) using spline ansatz functions by extending the approach from \cite{koellermeier2020spline} to the models in \cite{KowalskiTorrilhon2019}, including an analysis of hyperbolicity, a hyperbolic regularization for hyperbolic SSWME, and numerical test cases to assess accuracy and robustness. Similar to the spline-based models for rarefied gases \cite{koellermeier2020spline}, we show that spline basis functions increase the flexibility of moment models and yield accurate models for shallow flows. 

The rest of this paper is organized as follows: 
In \Fref{sec:2} we briefly present the free-surface reference system and explain the process of model reduction using moments. We derive constrained spline basis functions with zero-mean in \Fref{sec:constrainedsplines} and the new spline-based SSWME in \Fref{sec:SplineMomentSystem}, while relating the new models to existing ones before studying their hyperbolicity and performing a hyperbolic regularization with analytical proof of hyperbolicity for hierarchies of high-order SSWME. \Fref{sec:Experiments} shows extensive numerical results for a smooth wave and a fast wave scenario.
The paper ends with a conclusion and possibilities for future work.

\section{Derivation of Shallow Water Moment Equations}
\label{sec:2}
The derivation of reduced free-surface models in \cite{KowalskiTorrilhon2019} employs a reference system, obtained from the incompressible Euler equations by mapping onto a scaled vertical coordinate $\zeta$ and a subsequent Legendre basis expansion of the velocity profile in vertical direction, the so-called moment expansion. A Galerkin projection of the momentum equation onto test functions then leads to a closed set of equations, the so-called Shallow Water Moment Equations (SWME). In this section, we briefly recall those ingredients, but for a general basis function that will only be specified in the next section.

\subsection{Reference System}
\label{sec:ReferenceSystem}
Starting point are the incompressible Euler equations expressing the conservation of mass and momentum. For conciseness, only one horizontal direction is considered. The extension to two horizontal dimensions  is straightforward and given in \cite{KowalskiTorrilhon2019}. The equations, written in component form, read
\begin{align}
\label{eq:eMomentumHorizontal}
\partial_t u+\partial_x (u^2)+\partial_z(u w) &=
 - \frac{1}{\rho} \partial_x p+\frac{1}{\rho} \partial_x \sigma,\\
 \label{eq:eMomentumVertical}
 \partial_t w + \partial_x (u w) 
 +\partial_z (w^2) &=- \frac{1}{\rho} \partial_z p+\frac{1}{\rho} \partial_z \sigma - g,\\
 \label{eq:eMass}
             \partial_x u + \partial_z w &= 0,
\end{align}
where $u$ and $w$ denote the flow velocity in horizontal and vertical direction, respectively. Other variables are the pressure $p$ and the deviatoric stress tensor $\sigma$, depending on space $(x,z)$ and time $t$. In our setting, we assume the density $\rho$ to be constant and that the direction of gravitational acceleration $g$ is aligned with the $z$-axis.
We further assume a non-moving bottom $h_b$, i.e., $\partial_t h_b =0 $. Kinematic boundary conditions hold at the bottom and at the moving surface $h_s$:
\begin{align}
    \label{eq:bcwtop}
    &\text{at $z=h_s$:  }\partial_t h_s + u \partial_x h_s  = w,\\
    \label{eq:bcwbottom}
    &\text{at $z=h_b$:  }u\partial_x h_b  = w.
\end{align}
The difference $h_s-h_b$ corresponds to the fluid height and will be denoted by $h$ with $h>0$.

The pressure is assumed hydrostatic 
\begin{align}
    \label{eq:hydrostaticpressure} 
    p=( h_s-z)g,
\end{align}
and a Newtonian closure is chosen for the deviatoric stress tensor with kinematic viscosity $\nu$
\begin{align}
    \sigma_{xz}=\nu \rho\, \partial_z u .
\end{align}

Mapped onto $\zeta$-coordinates $\zeta = \frac{z-h_b}{h}$ \cite{ KowalskiTorrilhon2019,Phillips1956}, the system will take a different form: 
Due to the hydrostatic pressure \eqref{eq:hydrostaticpressure}, the momentum equation for the vertical velocity $w$ \eqref{eq:eMomentumVertical} decouples and $w$ can be restored using the divergence constraint \eqref{eq:eMass}.
The mapped system as derived in \cite{KowalskiTorrilhon2019} is called the reference system and reads 
\begin{align}
\label{eq:mappedmass}
\partial_{t}h+\partial_{x}\left(hu_{m}\right)+\partial_{y}\left(hv_{m}\right) & =0,\\
\label{eq:mappedmomentum}
\partial_t (h u + \frac{g}{2} h^2)
+ \partial_{\zeta}\left(h u \omega - \frac{\nu}{h} \partial_{\zeta} u\right)
& =g h \partial_x h_b.
\end{align}
The vertical coupling $\omega$ is defined as
\begin{align}
\omega & =\frac{1}{h}\int_{0}^{\zeta}\left(\int_{0}^{1}\partial_x( h u(\check{\zeta}))d\check{\zeta}-\partial_x( h u(\hat{\zeta}))\right)d\hat{\zeta}.
\label{eq:omega2-1}
\end{align}
The mapping also produces a new set of mapped boundary conditions
    \begin{align}
    \left.\partial_{\zeta}{u}\right|_{\zeta=1}=0\qquad
    \text{and}
    \qquad
    \left.\partial_{\zeta}{u}\right|_{\zeta=0}=\frac{h}{\lambda}\left.{u}\right|_{\zeta=0},
    \label{eq:mappedb}
    \end{align}
    where $\lambda$ stands for the slip length.

The reference system of equations \eqref{eq:mappedmass}-\eqref{eq:mappedmomentum} with boundary conditions \eqref{eq:mappedb} uses $\zeta$ instead of the previous vertical $z$.
In this formulation, the bottom is denoted by $\zeta=0$ and the free surface by $\zeta=1$.
These fixed boundaries are convenient for the application of a Galerkin projection method involving depth integration in the next subsection.

\subsection{Moment Expansion}
To reduce the dimension of the reference system \eqref{eq:mappedmass}-\eqref{eq:mappedmomentum}, a general ansatz for the horizontal velocity $u(x,t,\zeta)$ is used. This velocity will be represented by its vertical mean $u_m(x,t)=\int_0^1 u(x,t,\zeta)d\zeta$ plus an expansion around that mean
\begin{align}
    u(x,t,\zeta)=u_m(x,t)+ \sum_{i=1}^N s_i(x,t) \phi_i(\zeta).
\label{eq:momentansatz}
\end{align}
where the expansion uses suitable basis functions $\phi_i(\zeta)$ (yet to be defined) and basis coefficients denoted by $s_i$. 
Note that we use $s_i$ instead of $\alpha_i$ as introduced in \cite{KowalskiTorrilhon2019}, to distinguish our later spline model from the classical Shallow Water Moment Equations (SWME). These first models derived in \cite{KowalskiTorrilhon2019} used orthogonal Legendre basis functions, i.e., $\phi_i = \phi_i^{Leg}$ with $\phi_i^{Leg}$ the scaled and shifted Legendre polynomial of degree $i$, orthogonal on the integration domain $[0,1]$. 

However, the set of basis functions $\{\phi_i\}_{i=1,\ldots, N}$ could contain any functions~$\phi_i$ on $[0,1]$ as long as they are linearly independent. To ensure well-posedness and compatibility with the mean velocity $u_m$ the expansion part has to have zero-mean, denoted by
\begin{align}
\label{eq:zeroIntegralMean}
    \int_0^1 \sum_{i=1}^N s_i \phi_i(\zeta) d\zeta = 0 .
\end{align}
Enforcing \eqref{eq:zeroIntegralMean} by restricting the coefficients $s_i$ is briefly explained in \cite{Kauf2011,koellermeier2020spline}. Alternatively, choosing basis functions $\phi_i$ that separately fulfill the zero-mean condition, i.e. 
\begin{align}
    \int_{0}^1 \phi_i(\zeta)d\zeta=0 \textrm{ for all } i,
    \label{eq:zeroIntegralMeanSeparate}
\end{align}
is a stronger condition that will guarantee \eqref{eq:zeroIntegralMean} to hold. For the structure and sparsity of the resulting system of equations, it is beneficial if many
of the basis functions are orthogonal with respect to the $L_2$ scalar product $\langle f,g\rangle:=\int_{0}^1 f(\zeta) g(\zeta) \text{ d}\zeta$. Consequently, \cite{KowalskiTorrilhon2019} employs scaled versions of the first
$N$ Legendre polynomials, which are pairwise orthogonal and indeed fulfill the zero-mean condition. 

\subsection{Galerkin Projection and Moment System}
To derive a set of PDEs for the expansion coefficients, the ansatz \eqref{eq:momentansatz} is first inserted into the momentum equation of the reference system \eqref{eq:mappedmomentum}. This is equivalent to assuming that the vertical velocity profiles can be represented in the form of \eqref{eq:momentansatz}. 
Galerkin projections onto a set of suitable test functions then lead to a set of equations for the coefficients. The Galerkin projection is performed by integrating along the vertical axis a product of the equation \eqref{eq:mappedmomentum} with a test function, where the set of test functions is equal to the set of ansatz/basis functions $\{\phi_i\}_{i=1,\ldots, N}$ extended by the constant function $\phi_0 = 1$. The procedure is detailed in \cite{KowalskiTorrilhon2019}, where also the weak enforcement of the boundary conditions is described. Together with the mass balance \eqref{eq:mappedmass} this constitutes a closed system of equations. The full set of variables then consist of water height $h$, mean velocity $u_m$, and basis coefficients $s_i$, for $i=1,\ldots, N$, i.e., $N+2$ variables.

The resulting equation system of $N+2$ equations for the unknowns $(h,h u_m, h s_i)^T \in \mathbb{R}^{N+2}$ consists of three parts: one mass balance
\begin{equation}
    \partial_t h + \partial_x (h u_m) = 0,
    \label{eq:conteqn}
\end{equation}
which is followed by one averaged momentum equation
\begin{equation}
    \partial_t (h u_m)+ \partial_x \left( h u_m^2 + h \sum_{i,j=1}^N s_i s_j M_{ij} + \frac{g h^2}{2} \right)=gh \partial_x h_b-\frac{\nu}{\lambda} \left( u_m + \sum_{j=1}^N s_j V^0_j\right),
    \label{eq:momentumeqn}
\end{equation}
and finally $N$ Galerkin projections produce $N$ additional equations for the coefficients $s_i$
\begin{multline}
    \partial_t \left( h \sum_{j=1}^N s_j M_{ij} \right)
    + \partial_x \left( h \left( 2 u_m \sum_{j=1}^N M_{ij} s_j + \sum_{j,k=1}^N  A_{ijk} s_j s_k \right) \right) 
    - u_m \sum_{j=1}^N \partial_x (h s_j) M_{ij} \\
    + \sum_{j,k=1}^N B_{ijk} s_k \partial_x (h s_j) 
    =-\frac{\nu}{\lambda} V^0_i \left(u_m + \sum_{j=1}^N s_j V^0_j\right) - \frac{\nu}{h} \sum_{j=1}^N s_j C_{ij},
    \label{eq:momenteqn}
\end{multline}
where the matrix $M \in \mathbb{R}^{N \times N}$, the tensors $A,B \in \mathbb{R}^{N \times N \times N}$ and the vector $V^0\in \mathbb{R}^N$ encode properties of the chosen basis functions. Their definitions can be found in Subsection \ref{sec:Matrices} of the appendix. We note that $M, A, B, V^0$ can typically be computed a-priori analytically \cite{KoellermeierRominger2020}.
The entries $M_{ij}$ of the matrix $M$ are precisely the values of the $L_2$ scalar product $\langle \phi_i, \phi_j \rangle$ of the basis functions and therefore zero for $i\neq j$ if the basis is orthogonal, e.g., for the Legendre basis in \cite{KowalskiTorrilhon2019}.
An important difference of our more general approach in comparison to \cite{KowalskiTorrilhon2019} is
that the matrix $M$ does not have to be diagonal, leading to several coupled time derivative terms in \eqref{eq:momentumeqn} and \eqref{eq:momenteqn}. 

It is worth mentioning that $V_i^0 = \phi_i(0)$, i.e., the evaluation of the basis function at the bottom. This is why the term $u_b = u_m + \sum_{i=1}^N s_i V^0_i$ is the velocity profile \eqref{eq:momentansatz} evaluated at the bottom $\zeta = 0$. The right hand side then contains this term as $-\frac{\nu}{\lambda} u_b$ modeling the friction with the bottom. The other friction term $- \frac{\nu}{h} \sum_{j=1}^N s_j C_{ij}$ is modeling Newtonian friction and takes into account the derivatives of the basis functions via the definition of $C_{ij}$ in Subsection \ref{sec:Matrices}.

In matrix-vector notation for the variable vector $U=(h, hu_m, hs_1,\ldots, hs_N)^T \in \mathbb{R}^{N+2}$, the equations \eqref{eq:conteqn}, \eqref{eq:momentumeqn} and \eqref{eq:momenteqn} can be written as 
\begin{equation}\label{eq:SSWME_con}
    \begin{pmatrix}
        1 & &\\
         & 1 & \\
         & & M
    \end{pmatrix}\frac{\partial U}{\partial t} + A \frac{\partial U}{\partial x} = P
\end{equation} 
using the system matrix $A \in \mathbb{R}^{(N+2)\times(N+2)}$ defined by
\begin{equation}\label{eq:A}
    A =
    \left(\begin{array}{ccccc}
    0 & 1 & 0 & \cdots & 0 \\
    - u_m^2 + gh - \sum\limits_{j,k=1}^{N} M_{jk} s_j s_k & 2 u_m &  2 \sum\limits_{j=1}^{N} M_{1j} s_j & \cdots & 2 \sum\limits_{j=1}^{N} M_{Nj}s_j \\
    -2 u_m \sum\limits_{j=1}^{N} M_{1j} s_j - \sum\limits_{j,k=1}^{N} A_{1jk} s_j s_k & 2 \sum\limits_{j=1}^{N} M_{1j} s_j &  & & \\
    \vdots & \vdots & & \mathcal{A} &  \\
    -2 u_m \sum\limits_{j=1}^{N} M_{Nj} s_j - \sum\limits_{j,k=1}^{N} A_{Njk} s_j s_k & 2 \sum\limits_{j=1}^{N} M_{Nj} s_j &  & &
    \end{array}\right),
\end{equation}
where the lower right block matrix $\mathcal{A} \in \mathbb{R}^{N \times N}$ is defined by
\begin{equation}\label{eq:A_block}
    \mathcal{A}_{i,l} = \sum\limits_{j=1}^{N} \left( B_{ijl} + 2 A_{ijl} \right) s_j + u_m M_{il}
\end{equation}
and the right hand side friction vector $P$ is defined as
\begin{equation}\label{eq:P}
    P = \underbrace{\begin{pmatrix}
        0\\
        gh \partial_x h_b \\
        0 \\
        \vdots \\
        0
    \end{pmatrix}}_{P_0} - 
   \underbrace{\frac{\nu}{\lambda} (u_m+\sum_{j=1}^N s_j V^0_j)
    \begin{pmatrix}
        0\\
        1 \\
        V^0_1 \\
        \vdots \\
        V^0_N
    \end{pmatrix}}_{P_1}
    - \underbrace{\frac{\nu}{h}
    \begin{pmatrix}
        0\\
        0 \\
        \sum_{j=1}^N s_j C_{1j} \\
        \vdots \\
        \sum_{j=1}^N s_j C_{Nj}
    \end{pmatrix}}_{P_2}.
\end{equation}
Note that $P$ in \eqref{eq:P} contains three parts with different physical interpretations: (1) a bottom topography term $P_0$, (2) a bottom friction term $P_1$, and (3) an interior friction term $P_2$.  Alternatively, the system takes the form of a balance law 

\begin{equation}
    \label{eq:SSWME_con_flux}
    \begin{pmatrix}
        1 & &\\
         & 1 & \\
         & & M
    \end{pmatrix}\frac{\partial U}{\partial t} +  \frac{\partial F  (U)}{\partial x} =\tilde Q + P,
\end{equation}
with the friction term $P$ from \eqref{eq:P} and with the fluxes
\begin{equation}
    \label{eq:F}
    F(U)=\begin{pmatrix}
        h u_m \\
        \frac{g}{2}h^2+hu_m^2+h\sum\limits_{i,j=1}^N s_is_jM_{ij}\\
        h (2 u_m \sum\limits_{j=1}^N s_j M_{1j}+\sum\limits_{j,k=1}^N A_{1jk}s_j s_k) \\
        \vdots \\
        h (2 u_m \sum\limits_{j=1}^N s_j M_{Nj}+\sum\limits_{j,k=1}^N A_{Njk}s_j s_k)
    \end{pmatrix}
\end{equation}
and the nonconservative terms 
\begin{equation}
    \label{eq:Q}
    \tilde Q=\begin{pmatrix}
        0 \\
        0 \\
        u_m \sum\limits_{j=1}^NM_{1j} \frac{\partial(hs_j)}{\partial x}+\sum\limits_{j,k=1}^N B_{1jk} s_k \frac{\partial(hs_j)}{\partial x}\\
        \vdots\\
        u_m \sum\limits_{j=1}^NM_{Nj} \frac{\partial(hs_j)}{\partial x}+\sum\limits_{j,k=1}^N B_{Njk} s_k \frac{\partial(hs_j)}{\partial x}
    \end{pmatrix}.
\end{equation}
 
\section{Constrained spline basis functions}
\label{sec:constrainedsplines}
Global polynomial approximations require high regularity. Instead, approximations using splines, i.e., piecewise polynomials, are more robust for less regular data and therefore have a long history in engineering applications. This can be especially useful for shallow water profiles, which can include large derivatives or even layered structures, \cite{Fernandez-Nieto2016,garresdiaz2023}. Piecewise continuous spline functions have a polynomial degree $K \in \mathbb{N}$ and the interfaces of the polynomial pieces are defined on a grid $G=\left\{\zeta_0,\zeta_1, \dots, \zeta_{M-1}\right\}$ with nodes $\zeta_0=0 \leq \zeta_1 \leq \dots \leq \zeta_{M-1}=1$. At those nodes, each spline fulfills continuity conditions including its first $K-1$ derivatives. 

For an equidistant grid holds: The more grid points or the narrower the grid along the vertical axis and the higher the allowed degree $K$ of the polynomial pieces the better the approximation qualities of the spline basis. Each spline basis is therefore characterized by a number of grid points $M$ and an allowed degree $K$.

\begin{figure}
    \centering
    \includegraphics[width=.25\textwidth]{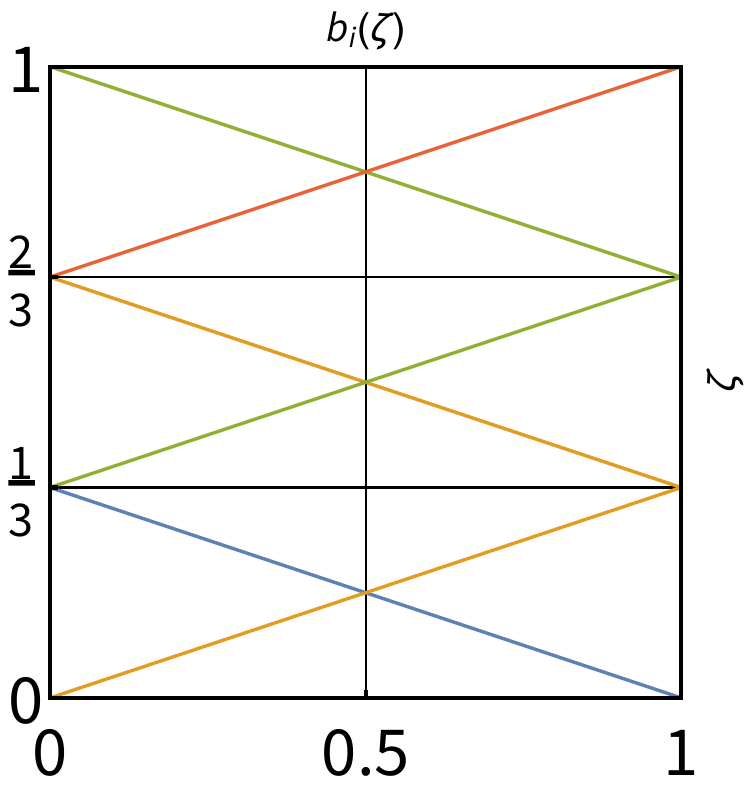}~~
    \includegraphics[width=.25\textwidth]{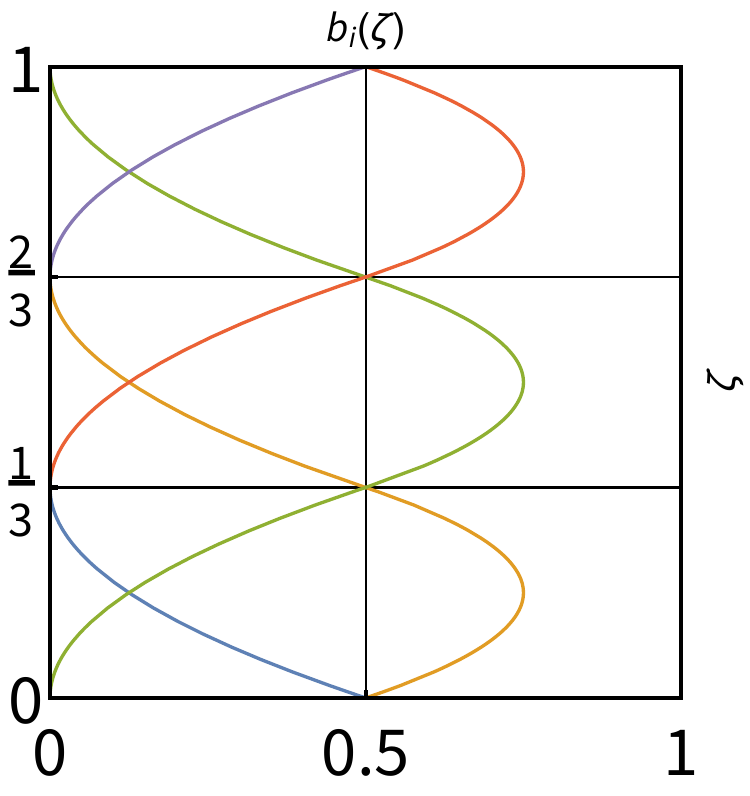}
    \caption{Linear B-splines (left) and quadratic B-splines (right) on an equidistant grid with $M=4$ points, i.e., $\Delta \zeta=\frac{1}{3}$. The grid allows for $M=4$ linear splines and $M+1=5$ quadratic splines.}
    \label{fig:basis}
\end{figure}

We thus model the velocity profile $u(x,t,\zeta)$ along $\zeta$ not as a global polynomial like \cite{KowalskiTorrilhon2019}, but as a piecewise polynomial. One could simply choose a grid and select the corresponding B-spline functions $\{b_i\}_{i=1,\ldots, M+K-1}$ as basis functions in \eqref{eq:momentansatz}. It is known that for a fixed degree~$K$ the B-splines form a basis for piecewise polynomial functions of degree $K$ on the predefined grid $G$ \cite{koellermeier2020spline}. Two typical sets of B-splines on a regular grid are depicted in Figure \ref{fig:basis}: linear and quadratic B-splines. Note that the number of quadratic B-splines on the same grid is one larger compared to linear splines. In general, a grid with $|G|=M$ grid points and B-splines of degree $K$ has an associated basis consisting of $M+K-1$ B-splines. This means $M$ linear B-splines ($K=1$) and $M+1$ quadratic B-splines ($K=2$).

However, this expansion using the standard B-splines does not fulfill the zero-mean condition \eqref{eq:zeroIntegralMeanSeparate},
which contradicts the assumption that $u_m$ is the vertical average of $u$ in the expansion \eqref{eq:momentansatz}. 
We therefore formulate new ansatz functions, the so-called constrained splines $\phi_i(\zeta)$, where  
\begin{equation}
\label{eq:zero-mean-phii}
    \int_0^1 \phi_i(\zeta) d \zeta = 0
\end{equation} is guaranteed by construction based on a linear combinations of two B-splines:
\begin{figure}
    \centering
    \includegraphics[width=.25\textwidth]{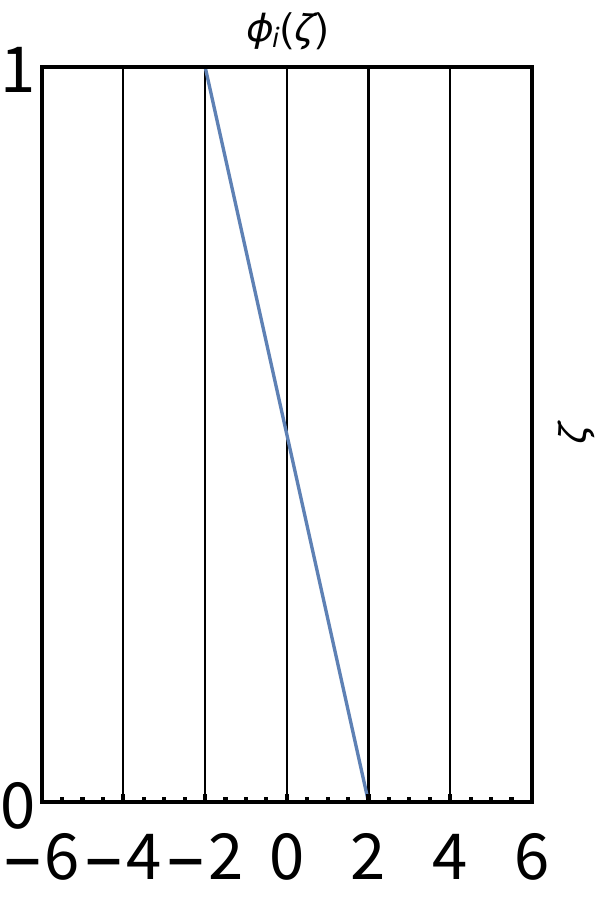}~~
    \includegraphics[width=.25\textwidth]{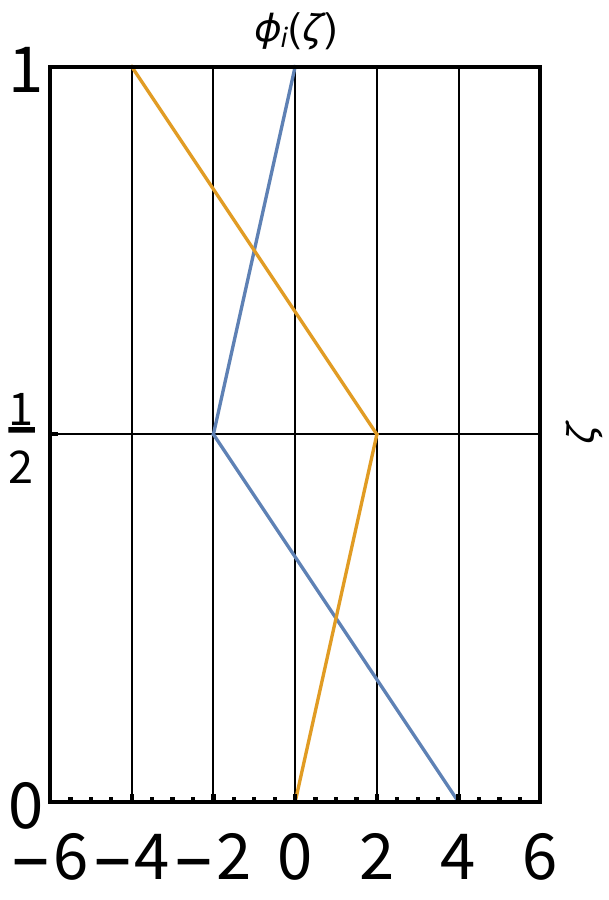}~~
    \includegraphics[width=.25\textwidth]{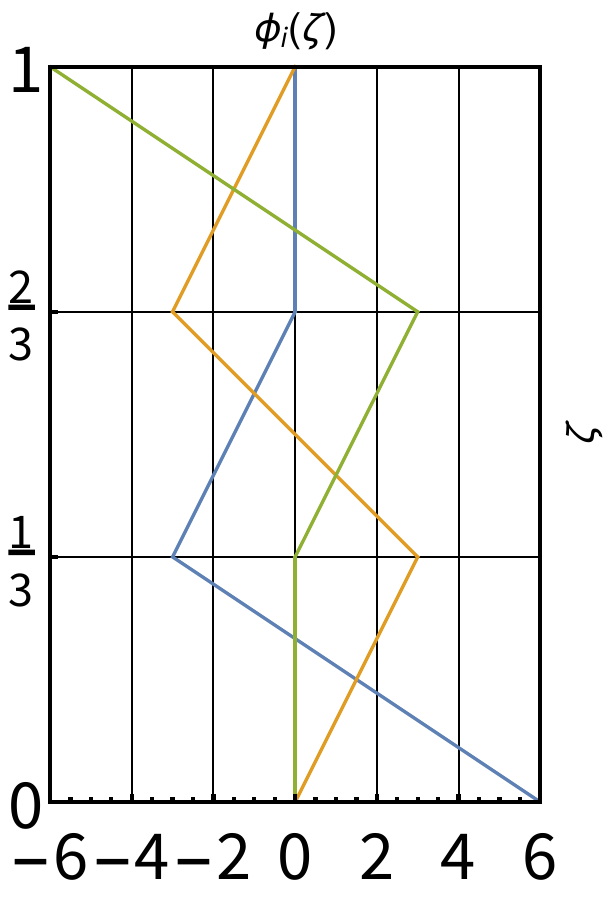} 
    \caption{Three constrained linear spline bases L1, L2, L3 satisfying \eqref{eq:zero-mean-phii} on equidistant grids: L1 uses $N=1$ spline (left); L2 uses $N=2$ splines (middle); L3 uses $N=3$ splines (right). Each linear spline basis uses $M=N+1$ grid points resulting in a grid width $\Delta \zeta = \frac{1}{N}$.}
    \label{fig:constrainedlinearbasis}
\end{figure}

\begin{figure}
    \centering
    \includegraphics[width=.25\textwidth]{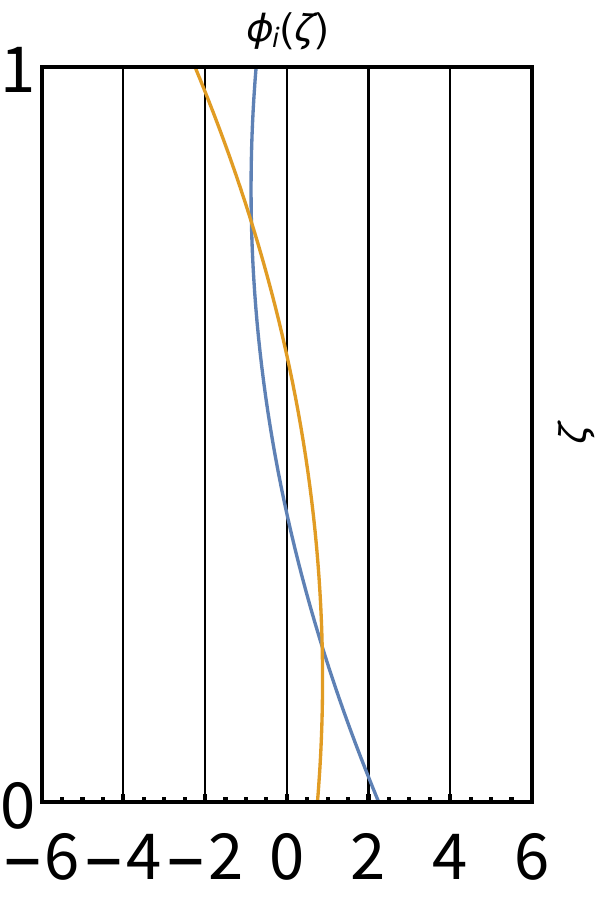}~~
    \includegraphics[width=.25\textwidth]{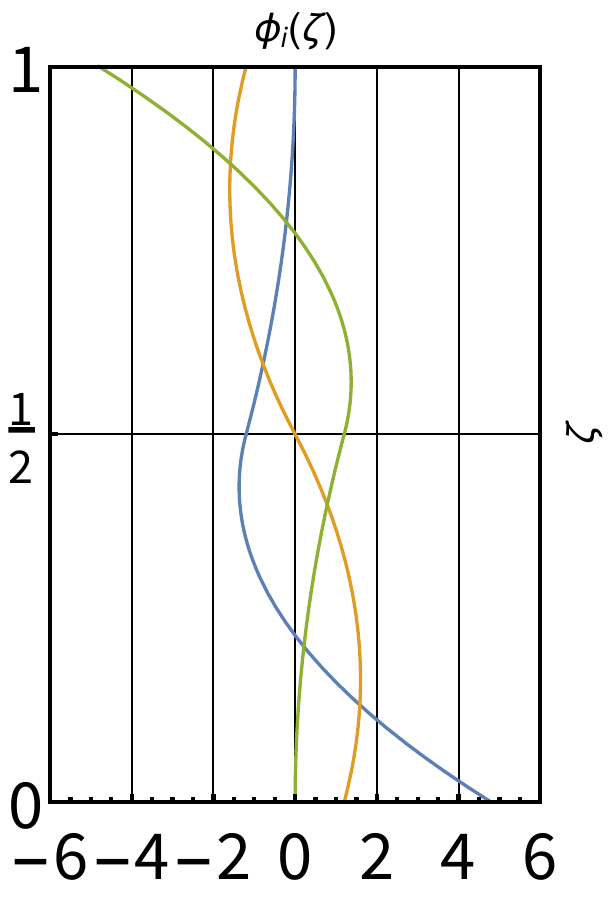}~~
    \includegraphics[width=.25\textwidth]{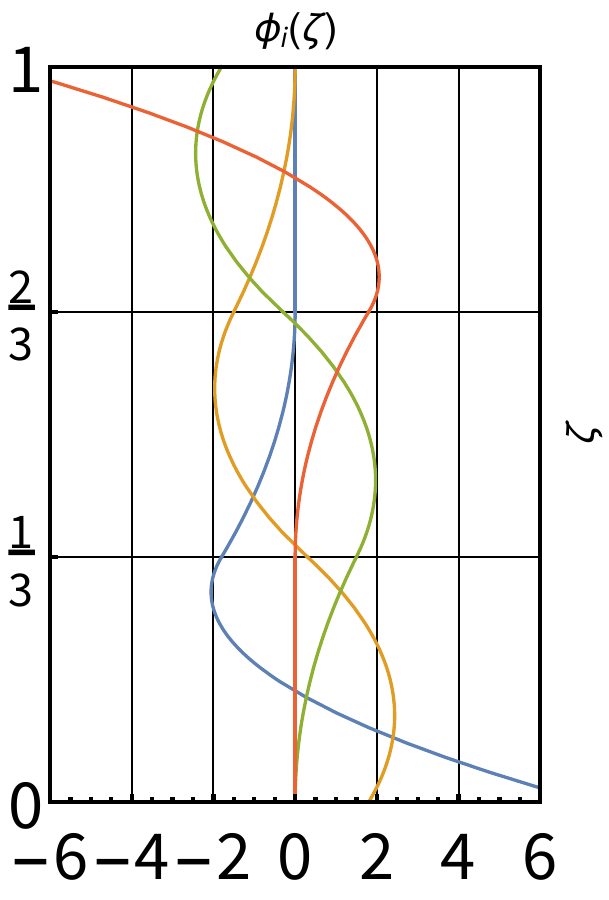}
    \caption{Three constrained quadratic spline bases Q2, Q3, Q4 satisfying \eqref{eq:zero-mean-phii} on equidistant grids: Q2 uses $N=2$ splines (left); Q3 uses $N=3$ splines (middle); Q4 uses $N=4$ splines (right). Each quadratic spline basis uses $M=N$ grid points resulting in a grid width $\Delta \zeta = \frac{1}{N-1}$.}
    \label{fig:constrainedquadraticbasis}
\end{figure}

Suppose a grid $G$ with B-splines $\{b_i\}_{i=1,\ldots, M+K-1}$ is given. The functions 
\begin{equation}
    \{\phi_i(\zeta):=b_i(\zeta) + \mu_i \, b_{i+1}(\zeta)\}, \quad {i=1,\ldots, M+K-2},
\end{equation} where $\mu_i$ is determined by the constraint \eqref{eq:zero-mean-phii}, form a new set of basis functions. The new ansatz functions are still linearly independent as 
the original B-splines are linearly independent. Due to \eqref{eq:zero-mean-phii} the constrained splines have the same regularity as the original B-splines.
By construction, the new basis $\{\phi_i\}_{i=1,\ldots, M+K-2}$ has one element less compared to the set of B-splines as it is missing functions offset by a non-zero mean. In particular, a basis containing $N$ linear constrained splines requires $M=N+1$ grid points which corresponds to a width of $\Delta \zeta=\frac{1}{N}$ on a regular grid and a basis containing $N$ quadratic constrained splines requires $M=N$ grid points implying a grid width of $\Delta \zeta=\frac{1}{N-1}$. Figures \ref{fig:constrainedlinearbasis} and \ref{fig:constrainedquadraticbasis} illustrate examples of the new constrained linear and constrained quadratic spline basis functions, respectively.

Below, we give explicit formulas for the cases $N=1, 2, 3$ for linear and quadratic constrained splines on an equidistant grid. 



\subsection{$N=1$ linear spline basis (L1)}
\label{sec:L1}
We first consider $N=1$ linear constrained basis function $\phi^{L1}_1$. The superscript here denotes the basis (linear, $N=1$) and the subscript denotes the index of the basis function (here only $i=1$). The single L1 basis function reads 
\begin{equation}
    \phi^{L1}_1(\zeta) = 2 - 4 \zeta,
    \label{eq:one_spline}
\end{equation}
which coincides with the linear Legendre polynomial $\phi^{Leg}_1$ used in \cite{KowalskiTorrilhon2019} up to a factor 2, i.e., $\phi^{L1}_1(\zeta) = 2 \phi^{Leg}_1(\zeta)$.
As a consequence, the span of possible velocity profiles in the resulting model
corresponds to that of the linear model in \cite{KowalskiTorrilhon2019}. L1 is depicted in Figure \ref{fig:constrainedlinearbasis} (left).


\subsection{$N=2$ linear splines basis (L2)}
\label{sec:L2}
Next, consider $N=2$ linear constrained basis functions $\phi^{L2}_1, \phi^{L2}_2$ forming the basis L2:
\begin{equation}
    \phi^{L2}_1(\zeta) =
    \begin{cases}
        4-12 \xi  & 0\leq \xi <\frac{1}{2} \\
        -4 + 4\xi & \frac{1}{2}\leq \xi \leq 1 
    \end{cases},\quad
    \phi^{L2}_2(\zeta) =
    \begin{cases}
        4 \xi  & 0\leq \xi <\frac{1}{2} \\
        8-12 \xi  & \frac{1}{2}\leq \xi \leq 1 
    \end{cases}
    \label{eq:two_linear_splines}
\end{equation}
The L2 basis is depicted in Figure \ref{fig:constrainedlinearbasis} (middle).
The asymmetric shape of the piecewise linear functions is typical for the constrained linear spline bases.

Note that the L1 basis is included in L2 as $\frac{1}{2}\left(\phi^{L2}_1 + \phi^{L2}_2\right) = \phi^{L1}_1$. In turn, this also means that the linear Legendre function $\phi_1^{Leg} = \frac{1}{4}\left(\phi^{L2}_1 + \phi^{L2}_2\right)$ is included in L2.

\subsection{$N=3$ linear splines basis (L3)}
\label{sec:L3}
We now consider $N=3$ linear constrained basis functions $\phi^{L3}_1, \phi^{L3}_2, \phi^{L3}_3$, as depicted in Figure \ref{fig:constrainedlinearbasis} (right). The L3 basis reads
\begin{equation}
    \phi^{L3}_1(\zeta) = 
    \begin{cases}
        6-27 \zeta  & 0\leq \zeta <\frac{1}{3} \\
        -6 + 9 \zeta & \frac{1}{3}\leq \zeta \leq \frac{2}{3} \\
        0 & \frac{2}{3}\leq \zeta \leq 1 
    \end{cases}, \quad
    \phi^{L3}_2(\zeta) = 
    \begin{cases}
        9 \zeta  & 0\leq \zeta <\frac{1}{3} \\
        9-18 \zeta  & \frac{1}{3}\leq \zeta <\frac{2}{3} \\
         -9 +9 \zeta & \frac{2}{3}\leq \zeta \leq 1 
    \end{cases}, \end{equation}
    \begin{equation*}
    \phi^{L3}_3(\zeta) = 
    \begin{cases}
        0  & 0\leq \zeta <\frac{1}{3} \\
         -3 +9 \zeta& \frac{1}{3}\leq \zeta <\frac{2}{3} \\
        21-27 \zeta  & \frac{2}{3}\leq \zeta \leq 1 
    \end{cases}.
    \label{eq:three_linear_splines}
\end{equation*}

Note that the $N=1$ linear basis function is included as $\frac{1}{3}\phi^{L3}_1 + \frac{5}{9}\phi^{L3}_2 + \frac{1}{3}\phi^{L3}_3 = \phi^{L1}_1$. 
Again, so is the linear Legendre function as $\phi_1^{Leg} = \frac{1}{6}\phi^{L3}_1 + \frac{5}{18}\phi^{L3}_2 + \frac{1}{6}\phi^{L3}_3.$

\subsection{$N=2$ quadratic splines basis (Q2)}
\label{sec:Q2}
We note that no quadratic constrained basis Q1 with only one basis function exists, since at least two degrees of freedom are necessary to represent all piecewise quadratic functions.

The Q2 basis consists of the $N=2$ constrained quadratic spline basis functions $\phi^{Q2}_1$ and $\phi^{Q2}_2$ as depicted in Figure \ref{fig:constrainedquadraticbasis} (left). The basis functions are
\begin{equation}
    \phi^{Q2}_1(\zeta) = 
    \frac{3}{4} \left(6 \zeta ^2-10 \zeta +3\right), \quad
    \phi^{Q2}_2(\zeta) =
    -\frac{3}{4} \left(6 \zeta ^2-2 \zeta -1\right).
    \label{eq:two_quadratic_splines}
\end{equation}

We first note that the $N=1$ linear constrained basis function is included as $\frac{2}{3}\left(\phi^{Q2}_1 + \phi^{Q2}_2\right) = \phi^{L1}_1$. This means that the linear Legendre polynomial is included in Q2 via $\phi^{Leg}_1 = \frac{1}{3}\left(\phi^{Q2}_1 + \phi^{Q2}_2\right)$.

Furthermore, the basis is equivalent to the second order Legendre basis $\phi^{Leg}_1,\phi^{Leg}_2$ used in \cite{KowalskiTorrilhon2019}
\begin{equation}
    \phi^{Leg}_1(\zeta) = 1-2 \zeta, \quad
    \phi^{Leg}_2(\zeta) = 6 \zeta^2-6 \zeta +1,
\end{equation}
because $\phi^{Leg}_1,\phi^{Leg}_2$ can be expressed as linear combinations of $\phi^{Q2}_1,\phi^{Q2}_2$ and vice versa via

\begin{equation}
    \left(\phi^{Leg}_1,\phi^{Leg}_2 \right)^T = T\left(\phi^{Q2}_1,\phi^{Q2}_2\right) \text{ where }
    T = \begin{pmatrix}
        \frac{1}{3} & \frac{1}{3} \\
        \frac{2}{3} & -\frac{2}{3}
    \end{pmatrix}
\end{equation}
and, respectively,
\begin{equation}
    \left(\phi^{Q2}_1,\phi^{Q2}_2\right)^T = T^{-1} \left(\phi^{Leg}_1,\phi^{Leg}_2\right)^T \text{ with }
    T^{-1}=\begin{pmatrix}
        \frac{3}{2} & \frac{3}{4} \\
        \frac{3}{2} & -\frac{3}{4} 
    \end{pmatrix}.
\end{equation}


\subsection{$N=3$ quadratic splines basis (Q3)}
\label{sec:Q3}
Lastly, the $N=3$ constrained quadratic spline basis functions $\phi^{Q3}_1,\phi^{Q3}_2,\phi^{Q3}_3$ result in the basis Q3 depicted in Figure \ref{fig:constrainedquadraticbasis} (middle) and read
\begin{equation}
   \phi^{Q3}_1(\zeta) = 
    \begin{cases}
 \frac{24}{5} \left(7 \zeta ^2-6 \zeta +1\right) & 0\leq \zeta < \frac{1}{2} \\
 -\frac{24}{5} (\zeta^2 -2 \zeta + 1) & \frac{1}{2}\leq \zeta \leq 1 

    \end{cases},~  
    \phi^{Q3}_2(\zeta) = 
    \begin{cases}
 -\frac{6}{5} \left(12 \zeta ^2-4 \zeta -1\right) & 0\leq \zeta
   <\frac{1}{2}  \\
       \frac{6}{5} \left(12 \zeta ^2-20 \zeta +7\right) & \frac{1}{2}\leq
   \zeta \leq 1 
    \end{cases}, 
    \label{eq:three_quadratic_splines}
    \end{equation}
\begin{equation*}
    \phi^{Q3}_3(\zeta) = 
        \begin{cases}
        \frac{24 }{5}\zeta ^2 & 0\leq \zeta <\frac{1}{2} \\
 -\frac{24}{5} \left(7 \zeta ^2-8 \zeta +2\right) & \frac{1}{2}\leq
   \zeta \leq 1  
    \end{cases}.
\end{equation*}

We note that the L1 constrained basis is included as $\frac{1}{4}\phi^{Q3}_1 + \frac{2}{3} \phi^{Q3}_2 + \frac{1}{4} \phi^{Q3}_3 = \phi^{L1}_1$. 
This means that the linear Legendre polynomial is included in Q3 using $\phi^{Leg}_1 = \frac{1}{8}\phi^{Q3}_1 + \frac{1}{3} \phi^{Q3}_2 + \frac{1}{8} \phi^{Q3}_3$. 
Also the Q2 basis is included as $\frac{11}{32}\phi^{Q3}_1 + \frac{1}{2} \phi^{Q3}_2 + \frac{1}{32} \phi^{Q3}_3 = \phi^{Q2}_1$ and $\frac{1}{32}\phi^{Q3}_1 + \frac{1}{2} \phi^{Q3}_2 + \frac{11}{32} \phi^{Q3}_3 = \phi^{Q2}_2$.

\begin{remark}
    As evident from the inclusions of 'smaller' spline bases in 'larger' spline bases shown above, the spline bases build a hierarchical set of bases that can be used to derive a hierarchy of corresponding models later on. The derivations above can be generalized as a larger spline model includes a smaller one if the grid points of the smaller grid are included in the grid points of larger grid and the degree of the splines is at least as large. A more detailed study of the hierarchical structure might pave a way to adaptive spline models that yield a local refinement in some portion of the domain $[0,1]$ by enriching the spline basis with the corresponding basis functions, similar to wavelets. However, this is out of the scope of this work and left for future studies.
\end{remark}

\section{Spline Shallow Water Moment Equations (SSWME)}
\label{sec:SplineMomentSystem}
After constructing the suitable constrained spline bases, they can be used for the ansatz \eqref{eq:momentansatz} to derive the specific moment models of the form \eqref{eq:momentumeqn}-\eqref{eq:momenteqn}. We exemplarily show the models corresponding to the bases L1, L2, L3 and Q2, Q3, as derived in the previous section. Models with larger bases can readily be derived in the same way.
As for the naming, we call the models Spline Shallow Water Moment Equations (SSWME) and append the basis abbreviation, e.g., starting with SSWME-L1 for linear splines with $N=1$. In contrast, we refer to the Legendre models of \cite{KowalskiTorrilhon2019} using the abbreviations SWME1, SWME2 etc.

We write each model derived from a basis B as a system of hyperbolic balance laws of the form  
\begin{equation}
\label{eq:balance_law}
    \partial_t U^B + \partial_x F^B(U^B) = Q^B(U^B) \partial_x U^B -P^B,
\end{equation}
where the superscript $(\cdot)^B$ denotes the basis $B$, e.g., $B \hat{=} L1$, indicating that all parts of the system change with the basis. 

Furthermore, \eqref{eq:balance_law} includes the variable vector $U^B=(h, hu_m, h s^B_1, \ldots, h s^B_N)^T \in \mathbb{R}^{N+2}$, the flux function $F^B \in \mathbb{R}^{N+2}$, the non-conservative matrix $Q^B \in \mathbb{R}^{(N+2)\times(N+2)}$ and the friction term $P^B \in \mathbb{R}^{N+2}$, as derived from \eqref{eq:SSWME_con_flux} (Note that the $P, P_1, P_2$ terms in this Section differ from those in Section \ref{sec:2}, as those are not multiplied by $M^{-1}$). The models thus differ in the set of variables, their conservative and non-conservative fluxes as well as the friction terms.
For easier readability, we will drop the superscripts B wherever possible, i.e., when it is clear which basis we refer to.

For simplicity, we further consider a flat bottom $h_b = const$, in \eqref{eq:momentumeqn} for the rest of this paper.

\subsection{$N=1$ linear spline model (SSWME-L1)}
\label{sec:SSWME-L1}
We consider the basis L1. In this simplest case of $N=1$ linear spline, there is only one additional coefficient $s_1$ of the linear constrained spline function $\phi^{L1}_1$ from \eqref{eq:one_spline}, plotted in Figure \ref{fig:constrainedlinearbasis} (left). The grid for the spline functions is trivial because the whole
domain [0,1] is one grid cell. This leads to a variable vector $U=(h,h u_m, h s_1)^T \in \mathbb{R}^3$. In Section \ref{sec:L1} we have shown that the basis 
corresponds to the first order Legendre basis of \cite{KowalskiTorrilhon2019} for $s_1 = \frac{1}{2} \alpha_1$. The same is true for the resulting system of equations. For completeness, we state the SSWME-L1 system using the notation of this paper:
\begin{equation}
\label{eq:SSWME-L1}
\partial_t \begin{pmatrix}
    h \\ h u_m\\ h s_1
\end{pmatrix}
+ \partial_x \begin{pmatrix}
 h u_m \\
 \frac{g}{2}h^2+h u_m^2+\frac{4 h s_1^2}{3} \\
 2 h s_1 u_m \\
\end{pmatrix}
    =
    \begin{pmatrix}
        0 \\ 0 \\ u_m \partial_x s_1
    \end{pmatrix}
    - \frac{\nu}{\lambda}\begin{pmatrix}
 0 \\
 u_m+2 s_1 \\
 \frac{3}{2} \left(u_m + 2 s_1\right) \\
\end{pmatrix}-\frac{\nu}{h}
\begin{pmatrix}
    0 \\
    0 \\
    12 s_1
\end{pmatrix}
\end{equation}
Note furthermore that the right hand side friction term includes the term $u_m+2 s_1$ corresponding to the velocity $u_b$ at the bottom. The remaining friction term models Newtonian friction inside the fluid due to non vanishing gradients of the velocity profile.

\subsection{$N=2$ linear splines model (SSWME-L2)}
\label{sec:SSWME-L2}
Consider now the basis L2 with $N=2$ linear splines. The vertical profile is approximated by a sum of two piecewise linear splines from \eqref{eq:two_linear_splines} as shown in Section \ref{sec:L2}, plotted in Figure \ref{fig:constrainedlinearbasis} (middle). The interior grid point is at $\zeta=\frac{1}{2}$. 
The SSWME-L2 system reads
\begin{equation}
\label{eq:SSWME-L2}
\partial_t 
\begin{pmatrix}
    h \\ h u_m \\ h s_1 \\ h s_2
\end{pmatrix}
+ \partial_x \begin{pmatrix}
h u_m \\
 \frac{g}{2}h^2+h u_m^2+\frac{8 h s_1^2}{3}+\frac{8 h s_2^2}{3} \\
 2 h s_1 u_m+\frac{3 h s_1^2}{2}+h s_2 s_1-\frac{h s_2^2}{2} \\
 2 h s_2 u_m+\frac{h s_1^2}{2}-h s_2 s_1-\frac{3 h s_2^2}{2}
\end{pmatrix}
= Q
\partial_x 
\begin{pmatrix}
    h \\ h u_m \\ h s_1 \\ h s_2
\end{pmatrix} - P_1 -P_2,
\end{equation}
where 
\begin{equation}
    Q=
\begin{pmatrix}
 0 & 0 & 0 & 0 \\
 0 & 0 & 0 & 0 \\
 0 & 0 & u_m+\frac{3 s_1}{4}+\frac{5 s_2}{4} & \frac{s_1}{4}+\frac{3 s_2}{4} \\
 0 & 0 & -\frac{3 s_1}{4}-\frac{s_2}{4} & u_m-\frac{5 s_1}{4}-\frac{3 s_2}{4} \\
\end{pmatrix}
\end{equation}
and 
\begin{equation}P_1=\frac{\nu}{\lambda}\begin{pmatrix}
 0 \\
 4 s_1+u_m \\
 \frac{3}{2} (4 s_1+u_m) \\
 0 \\
\end{pmatrix}, \quad
    P_2= \frac{\nu}{h}
\begin{pmatrix}
 0 \\
 0 \\
 {30 s_1}-{18 s_2} \\
 {30 s_2}-{18 s_1} \\
\end{pmatrix}
    .
\end{equation}
 As shown in Section \ref{sec:L2}, the L1 basis is included in the L2 basis. We can therefore also derive the SSWME-L1 model from the SSWME-L2 model by setting $s^{L2}_1=s^{L2}_2=\frac{1}{2}s^{L1}_1$ and adding the last two equations of 
\eqref{eq:SSWME-L2}, which results in \eqref{eq:SSWME-L1}.


\subsection{$N=3$ linear splines model (SSWME-L3)}
\label{sec:SSWME-L3}
The system for the L3 basis with $N=3$ linear splines is called SSWME-L3 and can be written in the same form of \eqref{eq:balance_law}. It uses $P = P_1 + P_2$ and the variables $U=(h,h u_m, h s_1, h s_2, h s_3)^T$, with
\begin{equation}
    F(U)=
\begin{pmatrix}
 h u_m \\
 \frac{g}{2}h^2+h u_m^2+4 h s_1^2+3 h s_2^2+4 h s_3^2-h s_1 s_2-h s_1 s_3-h s_2 s_3 \\
 2 h s_1 u_m+\frac{77 h s_1^2}{30}+\frac{26}{15} h s_2 s_1+\frac{1}{3} h s_3 s_1-\frac{h s_2^2}{2}-\frac{11 h s_3^2}{15}-\frac{1}{15} h s_2 s_3 \\
 2 h s_2 u_m+\frac{9 h s_1^2}{5}-\frac{6}{5} h s_2 s_1-\frac{9 h s_3^2}{5}+\frac{6}{5} h s_2 s_3 \\
 2 h s_3 u_m+\frac{11 h s_1^2}{15}+\frac{1}{15} h s_2 s_1-\frac{1}{3} h s_3 s_1+\frac{h s_2^2}{2}-\frac{77 h s_3^2}{30}-\frac{26}{15} h s_2 s_3 \\
\end{pmatrix},
\end{equation} 
\begin{equation}
    Q=
    \begin{pmatrix}
\mathbf{0_{2 \times 2}}&   & \mathbf{0_{2 \times 3}} &   \\
  & u_m+\frac{21 s_1}{20}+\frac{71 s_2}{40}-\frac{3 s_3}{40} & \frac{13 s_1}{40}+\frac{3\
 
  s_2}{4}-\frac{23 s_3}{40} & \frac{3 s_1}{40}-\frac{s_2}{40}-\frac{3 s_3}{10} \\
\mathbf{0_{3\times 2}} & -\frac{9 s_1}{10}-\frac{3 s_2}{40}+\frac{9 s_3}{40} & u_m-\frac{69 s_1}{40}+\frac{69\
 
  s_3}{40} & -\frac{9 s_1}{40}+\frac{3 s_2}{40}+\frac{9 s_3}{10} \\
  & \frac{3 s_1}{10}+\frac{s_2}{40}-\frac{3 s_3}{40} & \frac{23 s_1}{40}-\frac{3\
 
  s_2}{4}-\frac{13 s_3}{40} & u_m+\frac{3 s_1}{40}-\frac{71 s_2}{40}-\frac{21 s_3}{20} \\
\end{pmatrix},
\end{equation} 
and
\begin{equation}
P_1=\frac{\nu}{\lambda}\begin{pmatrix}
 0 \\
 6 s_1+u_m \\
 \frac{47}{30} (6 s_1+u_m) \\
 \frac{3}{10} (6 s_1+u_m) \\
 \frac{7}{30} (6 s_1+u_m) \\
\end{pmatrix}
     \text{ and }
     P_2=\frac{\nu}{h}\begin{pmatrix}
 0 \\
 0 \\
 \frac{324 s_1}{5 }-\frac{162 s_2}{5 }+\frac{54 s_3}{5 } \\
 -\frac{162 s_1}{5 }+\frac{216 s_2}{5 }-\frac{162 s_3}{5 } \\
 \frac{54 s_1}{5 }-\frac{162 s_2}{5 }+\frac{324 s_3}{5 } \\
\end{pmatrix}.
\end{equation}
Note that the SSWME-L1 system is recovered setting 
$s^{L3}_1=\frac{1}{3} s^{L1}_1$, $s^{L3}_2=\frac{5}{9} s^{L1}_1$, $s^{L3}_3=\frac{1}{3} s^{L1}_1$, and adding the last three equations of the SSWME-L3 which results in \eqref{eq:SSWME-L1}.

\subsection{$N=2$ quadratic splines model (SSWME-Q2)}
\label{sec:SSWME-Q2}
The smallest quadratic splines system uses Q2, i.e., $N=2$ quadratic splines. This SSWME-Q2 model is given by
\begin{equation}
\label{eq:SSWME-Q2}
    \partial_t \begin{pmatrix}
        h \\ h u_m \\ h s_1 \\ h s_2
    \end{pmatrix}
    + \partial_x
    \begin{pmatrix}
 h u_m \\
 \frac{g}{2}h^2+h u_m^2+\frac{69 h s_1^2}{80}+\frac{69 h s_2^2}{80}+\frac{51}{40} h s_1 s_2 \\
 2 h s_1 u_m+\frac{197 h s_1^2}{140}+\frac{25}{14} h s_2 s_1+\frac{113 h s_2^2}{140} \\
 2 h s_2 u_m-\frac{113}{140} h s_1^2-\frac{25}{14} h s_2 s_1-\frac{197 h s_2^2}{140} \\
\end{pmatrix} = Q \partial_x \begin{pmatrix}
    h \\ h u_m \\ h s_1 \\ h s_2
\end{pmatrix} -P,
\end{equation}
with the matrix for the non-conservative terms reads
\begin{equation}
    Q = \begin{pmatrix}
 0 & 0 & 0 & 0 \\
 0 & 0 & 0 & 0 \\
 0 & 0 & u_m+\frac{87 s_1}{56}+\frac{447 s_2}{280} & \frac{363 s_1}{280}+\frac{87 s_2}{56} \\
 0 & 0 & -\frac{87 s_1}{56}-\frac{363 s_2}{280} & u_m-\frac{447 s_1}{280}-\frac{87 s_2}{56} \\
\end{pmatrix}
\end{equation}
and the friction terms are
\begin{equation}P_1=\frac{\nu}{\lambda}\begin{pmatrix}
 0 \\
 \frac{9 s_1}{4}+\frac{3 s_2}{4}+u_m \\
 \frac{13}{3} \left(\frac{9 s_1}{4}+\frac{3 s_2}{4}+u_m\right) \\
 -\frac{7}{3} \left(\frac{9 s_1}{4}+\frac{3 s_2}{4}+u_m\right) \\
\end{pmatrix}, \quad
    P_2=\frac{\nu}{h}
    \begin{pmatrix}
 0 \\
 0 \\
 {36 s_1}{}-{24 s_2}{} \\
 {36 s_2}{}-{24 s_1}{} \\
\end{pmatrix}.
\end{equation}
As shown in Section \ref{sec:Q2}, the Q2 basis includes both the linear basis L1 and the quadratic Legendre basis ansatz from \cite{KowalskiTorrilhon2019}. The SSWME-Q2 model therefore also includes both the SSWME-L1 \eqref{eq:SSWME-L1} and the SWME2 model from \cite{KowalskiTorrilhon2019} using appropriate manipulations of the system (left out for brevity).

\subsection{$N=3$ quadratic splines model (SSWME-Q3)}
\label{sec:SSWME-Q3}
Using the Q3 basis with $N=3$ quadratic splines leads to the SSWME-Q3 model, which written in the form of \eqref{eq:balance_law}, uses $U=(h,hu_m,hs_1,hs_2,hs_3)^T$ and
\begin{equation}
    F(U) = 
    \begin{pmatrix}
 h u_m \\
 \frac{g}{2}h^2+h u_m^2+\frac{48 h s_1^2}{25}+\frac{204 h s_2^2}{125}+\frac{48 h s_3^2}{25}+\frac{156}{125} h s_1 s_2-\frac{96}{125} h s_1\
 
  s_3+\frac{156}{125} h s_2 s_3 \\
 2 h s_1 u_m+\frac{2161 h s_1^2}{875}+\frac{8987 h s_2 s_1}{3500}+\frac{2}{5} h s_3 s_1+\frac{17 h s_2^2}{70}-\frac{611 h s_3^2}{875}-\frac{2437 h s_2\
 
  s_3}{3500} \\
 2 h s_2 u_m+\frac{96 h s_1^2}{875}-\frac{492}{875} h s_2 s_1-\frac{96 h s_3^2}{875}+\frac{492}{875} h s_2 s_3 \\
 2 h s_3 u_m+\frac{611 h s_1^2}{875}+\frac{2437 h s_2 s_1}{3500}-\frac{2}{5} h s_3 s_1-\frac{17 h s_2^2}{70}-\frac{2161 h s_3^2}{875}-\frac{8987 h s_2\
 
  s_3}{3500} \\
\end{pmatrix},
\end{equation}
\begin{equation}
    Q_{in}=\begin{pmatrix}

   u_m+\frac{687 s_1}{350}+\frac{13191 s_2}{7000}-\frac{27 s_3}{250} & \frac{1263 s_1}{875}+\frac{279 s_2}{280}-\frac{963 s_3}{875} & \frac{77 s_1}{250}-\frac{2041 s_2}{7000}-\frac{377 s_3}{350} \\
 -\frac{304 s_1}{175}-\frac{578 s_2}{875}+\frac{64 s_3}{125} & u_m-\frac{1782 s_1}{875}+\frac{1782 s_3}{875} & -\frac{64 s_1}{125}+\frac{578 s_2}{875}+\frac{304 s_3}{175} \\
   \frac{377 s_1}{350}+\frac{2041 s_2}{7000}-\frac{77 s_3}{250} & \frac{963 s_1}{875}-\frac{279 s_2}{280}-\frac{1263 s_3}{875} & u_m+\frac{27 s_1}{250}-\frac{13191 s_2}{7000}-\frac{687 s_3}{350} \\
\end{pmatrix},
\end{equation} so that the non-conservative term is written as
\begin{equation}
    Q=\begin{pmatrix}
        \mathbf{0_{2\times2}} &  \mathbf{0_{2\times3}} \\
        \mathbf{0_{3\times2}} & Q_{in}
    \end{pmatrix} 
    \end{equation}
    and $P=P_1+P_2$, where
    \begin{equation}    
    P_1=\frac{\nu}{\lambda}\begin{pmatrix}
 0 \\
 u_m+\frac{24 s_1}{5}+\frac{6 s_2}{5} \\
 \frac{23}{8} \left(u_m+\frac{24 s_1}{5}+\frac{6 s_2}{5}\right) \\
 -\frac{2}{3} \left(u_m+\frac{24 s_1}{5}+\frac{6 s_2}{5}\right) \\
 \frac{19}{24} \left(u_m+\frac{24 s_1}{5}+\frac{6 s_2}{5}\right) \\
\end{pmatrix}, \quad
    P_2=\frac{\nu}{h}\begin{pmatrix}
 0 \\
 0 \\
 \frac{1274 s_1}{15 }-\frac{96 s_2}{5 }+\frac{374 s_3}{15 } \\
 -\frac{736 s_1}{15 }+\frac{144 s_2}{5 }-\frac{736 s_3}{15 } \\
 \frac{374 s_1}{15 }-\frac{96 s_2}{5 }+\frac{1274 s_3}{15 } \\
\end{pmatrix}.
\end{equation}
Again, suitable manipulations reduce the SSWME-Q3 system to the SSWME-L1 or SSWME-Q2.

\subsection{Propagation speeds and hyperbolicity}
A main benefit of the analytical spline models is that they allow for in-depth analysis and identification of terms in the model representing different physical effects, as was already exemplified for the right-hand side friction term. In addition, the closed form of the equations allows to derive propagation speeds of the models, which are obtained by first rewriting the system \eqref{eq:balance_law} in the form 
\begin{equation}
    \partial_t U + A_{sys} \, \partial_x U = P,
\end{equation}
by means of computing $A_{sys}=\frac{\partial F}{\partial U}-Q$. Subsequently, the eigenvalues of $A_{sys}$ yield the propagation speeds, relevant for transport of surface waves. If a full set of real eigenvalues exists, the system is called \emph{hyperbolic}, which is a desirable property as it is associated with stability of numerical methods and allows for physical interpretation of the propagation speeds \cite{Fan2016,Koellermeier2017d,KoellermeierRominger2020}.  

Since the SSWME-L1 model is equivalent to the first order Legendre model from \cite{KowalskiTorrilhon2019}, which is hyperbolic with eigenvalues $\{u_m, u_m+\sqrt{gh+\alpha_1^2},u_m-\sqrt{gh+\alpha_1^2}\}$, where $\alpha_1$ is the first moment in the Legendre model, we directly consider the SSWME-L2 model \eqref{eq:SSWME-L2}. Its system matrix reads
\begin{equation}
\label{eq:ASysSSWME-L2}
       A^{L2}_{sys} = \begin{pmatrix}
 0 & 1 & 0 & 0 \\
 gh-u_m^2-\frac{8 s_1^2}{3}-\frac{8 s_2^2}{3} & 2 u_m & \frac{16\
 
  s_1}{3} & \frac{16 s_2}{3} \\

 -\frac{3 s_1^2}{2}-s_1 s_2-2 s_1 u_m+\frac{s_2^2}{2}

  & 2 s_1 & u_m+\frac{9 s_1}{4}-\frac{s_2}{4} &\
 
  \frac{3}{4} s_1-\frac{7}{4} s_2 \\
 -2 s_2 u_m-\frac{s_1^2}{2}+s_2 s_1+\frac{3 s_2^2}{2} & 2 s_2 &\
 
  \frac{7}{4} s_1-\frac{3}{4} s_2 &  
 
  u_m+\frac{s_1}{4}-\frac{9}{4} s_2
\end{pmatrix},
  \end{equation}
and the eigenvalues take the form 
\begin{equation}
\label{eq:evform}
    \lambda = u_m \pm c \sqrt{gh},
\end{equation}
where $c$ is the root of a fourth order polynomial that depends on the scaled moments $\overline{s_1}=\frac{1}{\sqrt{gh}}s_1$ and $\overline{s_2}=\frac{1}{\sqrt{gh}}s_2$.
In this case the polynomial is
\begin{align}
    \nonumber
    &4 c^4+c^3 \left(10 \overline{s_2}-10 \overline{s_1}\right)+c^2\left(-35 \overline{s_1}^2-6 \overline{s_2}\overline{s_1}-35 \overline{s_2}^2-4\right)\\
    \label{eq:char_poly}
    &+c \left(16 \overline{s_1}^3-144 \overline{s_2} \overline{s_1}^2+144 \overline{s_2}^2 \overline{s_1}+10\overline{s_1}-16 \overline{s_2}^3-10 \overline{s_2}\right)\\
    \nonumber
    &-8 \overline{s_1}^4+48 \overline{s_2} \overline{s_1}^3+112 \overline{s_2}^2 \overline{s_1}^2+3\overline{s_1}^2+48 \overline{s_2}^3 \overline{s_1}+6 \overline{s_2} \overline{s_1}-8 \overline{s_2}^4+3 \overline{s_2}^2
\end{align}
The form \eqref{eq:evform} and \eqref{eq:char_poly} of the eigenvalues implies that the hyperbolicity, i.e., whether all eigenvalues are real-valued, does not depend on the average velocity $u_m$ and only implicitly on the height $h$ by means of the scaling of the moments $\overline{s_i}$. This is a feature that all SSWME moment models share, similar to the Legendre models from \cite{KowalskiTorrilhon2019}. 

For small values of the coefficients $s_{1}$, $s_2$ \eqref{eq:char_poly} has real roots rendering the system hyperbolic. This can be seen because of two reasons: (1) in the limit of $s_1=s_2=0$, the polynomial \eqref{eq:char_poly} simplifies to $c^4-c^2$ with real roots $c_{1,2}=0$ and $c_{3,4}=\pm 1$, and (2) the roots depend continuously on the coefficients of the polynomial. 
However, due to the dependence of the polynomial on $\overline{s_1}$ and $\overline{s_2}$, large values of these coefficients can lead to a loss of hyperbolicity, similar to the Legendre basis models from \cite{KowalskiTorrilhon2019}. 

In \cite{KoellermeierRominger2020}, a hyperbolic regularization was suggested that is based on a linearization around linear velocity profiles. For the Legendre models, this effectively means setting higher coefficients (which correspond to quadratic, cubic, ... portions of the velocity profile) to zero. This was shown to be hyperbolic and extensively studied and applied thereafter \cite{garres2021,Koellermeier2020Eq, VerbiestKoellermeier}. 

We adopt a similar idea for the spline bases here. For the SSWME-L2 model, according to \Fref{sec:SSWME-L2} we restrict the coefficients to $s^{L2}_1=s^{L2}_2= \frac{s^{L1}_1}{2}=: \frac{\alpha_1}{4}$, following the notation of $\alpha_1$ being the linear basis coefficient in \cite{KoellermeierRominger2020}. This effectively restricts the velocity profile to a linear function. For the system matrix, this results in 
\begin{equation}
\label{eq:ASysHSSWME-L2}
       A^{L2}_{sys,H} = 
\begin{pmatrix}
 0 & 1 & 0 & 0 \\
 -\frac{\alpha_1 ^2}{3}+gh-u_m^2 & 2 u_m & \frac{4 \alpha_1 }{3} & \frac{4 \alpha_1 }{3} \\
 -\frac{\alpha_1 ^2}{8}-\frac{\alpha_1  u_m}{2} & \frac{\alpha_1}{2}  & \frac{\alpha_1}{2}+u_m & -\frac{\alpha_1 }{4} \\
 \frac{\alpha_1 ^2}{8}-\frac{\alpha_1  u_m}{2} & \frac{\alpha_1}{2}  & \frac{\alpha_1 }{4} & u_m-\frac{\alpha_1}{2}  \\
\end{pmatrix},
  \end{equation}
where $\alpha_1$ is the linear coefficient of the best possible Legendre expansion fit to the spline basis (see next subsection for details). 

The eigenvalues of $A^{L2}_{sys,H}$ are of the form $\lambda = u_m \pm c \sqrt{gh}$ requiring roots of a polynomial depending on $\alpha_1=\overline{\alpha_1}\sqrt{gh}$ given by
\begin{equation}
\label{eq:charpolyH}
16c^4-c^2(19 \overline{\alpha_1} ^2 c^2-16)+
    {3} \overline{\alpha_1} ^4+\overline{\alpha_1} ^2 =0.
  \end{equation}
The roots are given by $\{\pm \frac{\sqrt{3} }{4} \overline{\alpha_1} , \pm \sqrt{\overline{\alpha_1} ^2+1}\}$, leading to a set of real eigenvalues $\{u_m \pm \frac{\sqrt{3}}{4} \alpha_1,u_m \pm\sqrt{gh+\alpha_1^2}$.

\begin{figure}
    \centering
    \includegraphics[width=0.45\linewidth]{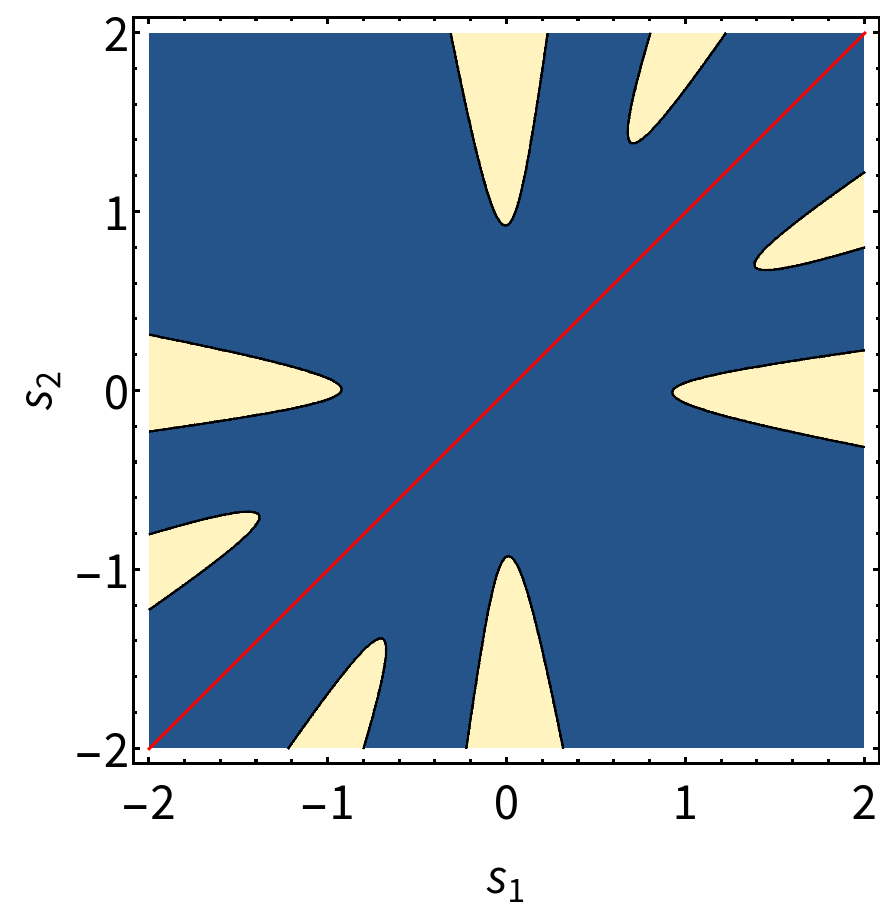}
    \includegraphics[width=0.45\linewidth]{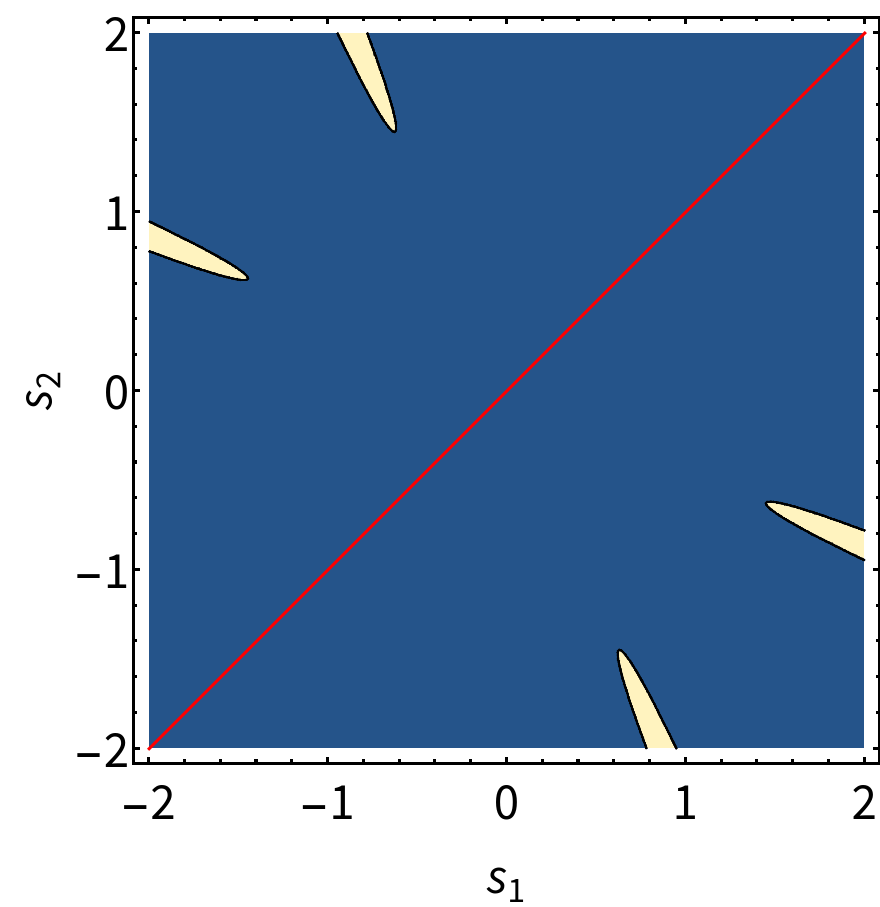}
    \caption{Hyperbolicity regions for SSWME-L2 (left) and SSWME-Q2 (right) depending on coefficients $s_1$ and $s_2$. The blue region corresponds to a locally hyperbolic system whereas in the yellow region imaginary eigenvalues of the flux Jacobian occur. 
    Enforcing a linear velocity profile restricts the coefficients to the red lines (here equivalent to $s_1=s_2$), along which both systems are hyperbolic, resulting in the HSSWME-L2 and HSSWME-Q2 models, respectively.}
    \label{fig:hyperbolicityneq2}
\end{figure}

The linearization of the system matrix around linear velocity profiles therefore creates a hyperbolic system from the SSWME-L2 model, which we denote by HSSWME-L2. The hyperbolicity is visualized in Figure \ref{fig:hyperbolicityneq2}. The left plot corresponds to the stated SSWME-L2 model, while the right plot shows the same approach for the SSWME-Q2 model. Yellow areas indicate loss of hyperbolicity due to the existence of roots with a non-zero imaginary part. 
There are extended regions of hyperbolicity in both models, with the SSWME-Q2 model notably showing an increased region of hyperbolicity around the origin.

The red lines forms a globally hyperbolic subspace, on the left by enforcing $s^{L2}_1=s^{L2}_2= \frac{s^{L1}_1}{2}=: \frac{\alpha_1}{4}$ for SSWME-L2 in the left plot.
For SSWME-Q2 the hyperbolic restriction that leads to a linear profile shape is $s^{Q2}_1=s^{Q2}_2=: \frac{\alpha_1}{3}$. This changes the system matrix from 
\begin{equation}
\label{eq:ASysSSWME-Q2}
        A^{Q2}_{sys}=\begin{pmatrix}
 0 & 1 & 0 & 0 \\
 gh-u_m^2-\frac{69 s_1^2}{80}-\frac{51 s_2 s_1}{40}-\frac{69 s_2^2}{80} & 2 u_m & \frac{69 s_1}{40}+\frac{51 s_2}{40} & \frac{51 s_1}{40}+\frac{69 s_2}{40} \\
 -2 s_1 u_m-\frac{197 s_1^2}{140}-\frac{25 s_2 s_1}{14}-\frac{113 s_2^2}{140} & 2 s_1 & u_m+\frac{353 s_1}{280}+\frac{53 s_2}{280} & \frac{137 s_1}{280}+\frac{17 s_2}{280} \\
 -2 s_2 u_m+\frac{113 s_1^2}{140}+\frac{25 s_2 s_1}{14}+\frac{197 s_2^2}{140} & 2 s_2 & -\frac{17 s_1}{280}-\frac{137 s_2}{280} & u_m-\frac{53 s_1}{280}-\frac{353 s_2}{280} \\
\end{pmatrix}
    \end{equation}
to 
\begin{equation}
\label{eq:ASysHSSWME-Q2}
    A^{Q2}_{sys,H}=\begin{pmatrix}
 0 & 1 & 0 & 0 \\
 -\frac{4 \alpha_1 ^2}{3}+gh-u_m^2 & 2 u_m & 2 \alpha_1  & 2 \alpha_1  \\
 -\frac{16 \alpha_1 ^2}{9}-\frac{4 \alpha_1  u_m}{3} & \frac{4 \alpha_1 }{3} & \frac{29 \alpha_1 }{30}+u_m & \frac{11 \alpha_1 }{30} \\
 \frac{16 \alpha_1 ^2}{9}-\frac{4 \alpha_1  u_m}{3} & \frac{4 \alpha_1 }{3} & -\frac{11 \alpha_1 }{30} & u_m-\frac{29 \alpha_1 }{30} \\
\end{pmatrix}
\end{equation} 
with eigenvalues $\{u_m \pm \frac{2 \alpha_1}{\sqrt{5}},u_m\pm\sqrt{gh+4 \alpha_1^2} \}$.
 In the right plot of Figure \ref{fig:hyperbolicityneq2} the restriction for the Q2 model is indicated by a red line, resulting in a hyperbolic model denoted as HSSWME-Q2.
 This shows that the linearization around linear velocity profiles is effective at restoring hyperbolicity in the resulting HSSWME-L2 and HSSWME-Q2 models.

\begin{figure}
    \centering
    \includegraphics[width=0.45\linewidth]{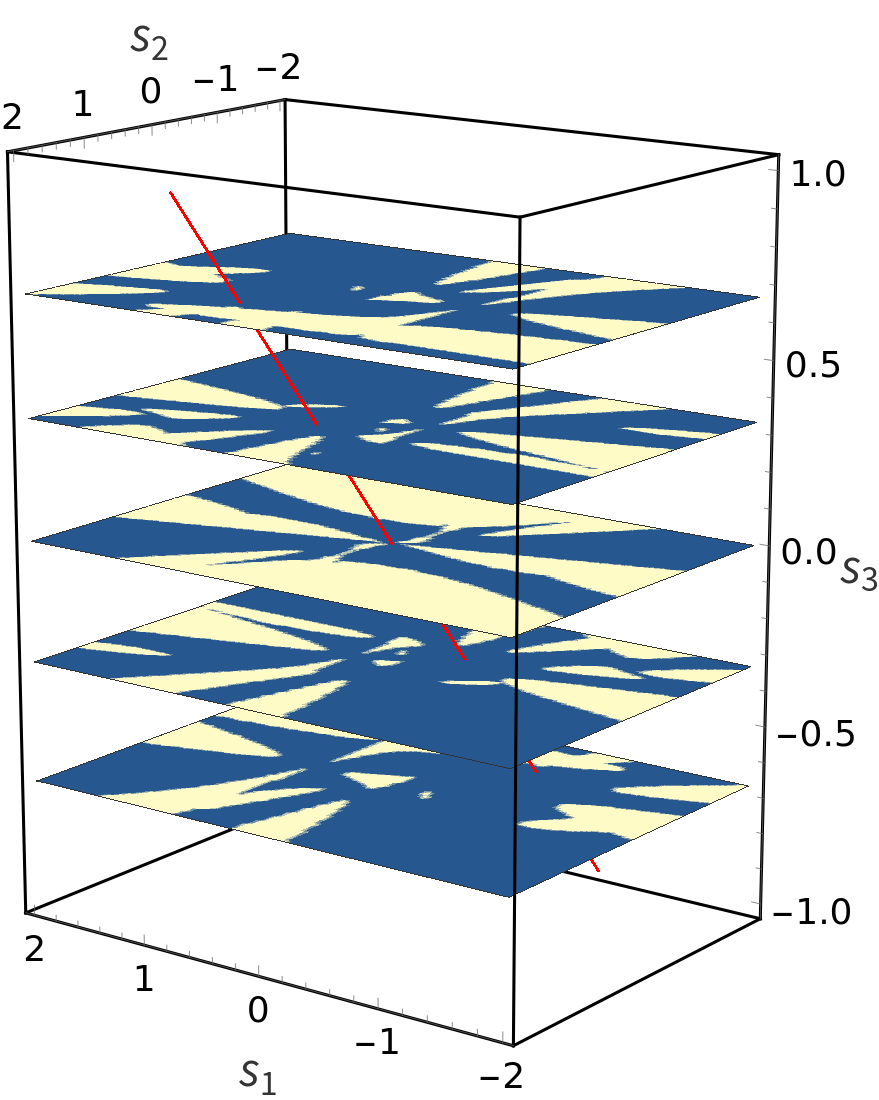}
    \includegraphics[width=0.45\linewidth]{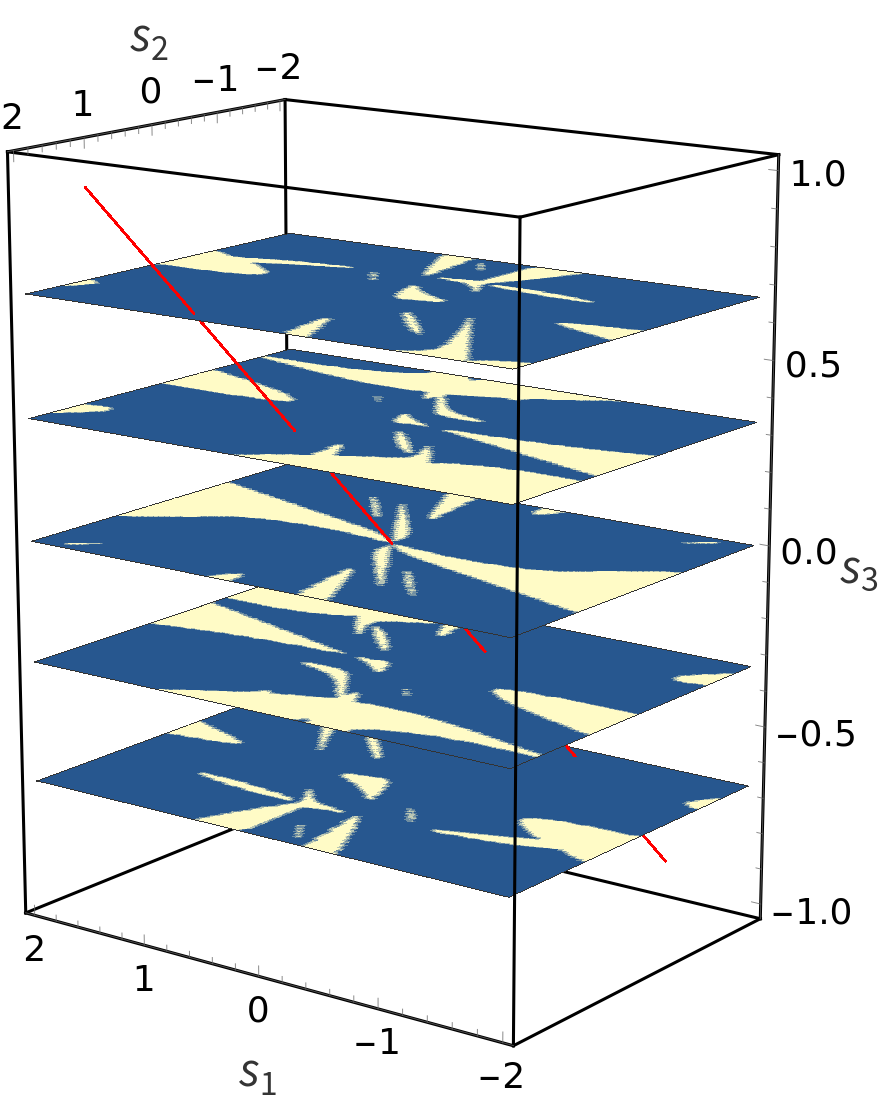}
    \caption{Hyperbolicity regions for SSWME-L3 (left) and SSWME-Q3 (right) depending on coefficients $s_1$, $s_2$, $s_3$. The blue region corresponds to a locally hyperbolic system whereas in the yellow region imaginary eigenvalues of the flux Jacobian occur. 
    Again, enforcing a linear velocity profile (here $s^{L2}_1=s^{L2}_2= \frac{s^{L1}_1}{2}=: \frac{\alpha_1}{4}$ (left) for L2 and $s^{Q2}_1=s^{Q2}_2=: \frac{\alpha_1}{3}$ (right) for Q2) restricts the coefficients to the red lines indicating a hyperbolic parameter subspace along which both systems are hyperbolic, resulting in the HSSWME-L3 and HSSWME-Q3 models, respectively.}
    \label{fig:hyperbolicityneq3}
\end{figure}

As for models with more basis functions, we consider the SSWME-L3 and SSWME-Q3 models. The $N=3$ basis coefficients render a visualization more difficult. We therefore show slices of the parameter domain in \Fref{fig:hyperbolicityneq3} (left) for SSWME-L3 and in \Fref{fig:hyperbolicityneq3} (right) for SSWME-Q3, respectively. Notably, both models with $N=3$ basis functions loose hyperbolicity for infinitely small deviations around the origin. This is similar to hyperbolicity loss in multi-dimensional moment models of rarefied gases \cite{Cai2014a} and motivates a further regularization of the systems in similar ways.
We therefore regularize the system matrix around linear velocity profiles by restricting the coefficients correspondingly. The details are given in Section \ref{sec:systems} of the appendix. This regularization restricts the system matrix evaluation to the red lines visible in \Fref{fig:hyperbolicityneq3}. Along these lines, the corresponding models traverse only the hyperbolic region. The regularized models called HSSWME-L3 and HSSWME-Q3 are thus globally hyperbolic.

\begin{remark}
    We denote that the hyperbolic regularization adopted here is not the only possibility of obtaining a hyperbolic model. Similar to the case of rarefied gas dynamics, where different regularizations exist in the literature \cite{Cai2014a,Cai2013b,Fan2016,Koellermeier2015,Koellermeier2014}, also Shallow Water Moment Models know different possible regularizations like the Shallow Water Linearized Moment Equations (SWLME) assuming small values of all coefficients only in the higher order equations \cite{Koellermeier2020Eq,Koellermeier2020steady}. A detailed investigation of other existing or potentially new hyperbolic regularization techniques for the SSWME of this paper is left for future work. The same is true for the derivation of corresponding energy/entropy equations, see, e.g., \cite{Gassner2016}.
\end{remark}

\subsection{Hyperbolicity proof for linear velocity profiles}
To prove hyperbolicity analytically for a hierarchy of models with arbitrary $N$, we will use the hyperbolicity result from \cite{KoellermeierRominger2020}, where it was shown that the system matrix is hyperbolic for a Legendre ansatz if the velocity profile is only linear.

We focus on a hierarchy of models with arbitrary $N$ where the spline basis consists of $N$ splines with piecewise degrees $N$. Examples from this hierarchy are the SSWME-L1 and the SSWME-Q2. The spline basis is then a linear transformation of the Legendre basis of degree $N$. 

The outline of the proof includes the following three steps: (1) We first transform the system matrix $A$ \eqref{eq:A} from the spline basis to the Legendre basis; (2) Next we perform the regularisation by projection onto linear velocity profiles; (3) Lastly, we transform back to the spline basis. As a similarity transformation of the hyperbolic Legendre model matrix, the regularized spline model matrix is then also hyperbolic.

For conciseness, we neglect the form of the friction term $P$ \eqref{eq:P}, as it is not relevant for the hyperbolicity of the model. A transformation of the friction term can be computed easily, since the friction term (for the model in this paper) is a linear function of average velocity $u_m$ and the moments $s_i$.

We denote with a superscript $S$ the SSWME system with spline basis functions $\phi^S_i$ similar as in \eqref{eq:SSWME_con} using variables
$(h, hu_m, hs)=(h, hu_m, hs_1,\ldots, hs_N) \in \mathbb{R}^{N+2}$ as 
\begin{equation}\label{eq:SSWME}
    \begin{pmatrix}
    1 & & \\ & 1 & \\ & & M^S
\end{pmatrix}\frac{\partial}{\partial t}\begin{pmatrix}
    h \\ h u_m \\ h s
\end{pmatrix} + A^S \frac{\partial}{\partial x}\begin{pmatrix}
    h \\ h u_m \\ h s 
\end{pmatrix} = 0,
\end{equation} 
while we denote with a superscript $L$ the SWME system from \cite{KowalskiTorrilhon2019} with Legendre basis functions $\phi^{Leg}_i$ using variables
$(h, hu_m, h \alpha)=(h, hu_m, h \alpha_1,\ldots, h \alpha_N) \in \mathbb{R}^{N+2}$ as 
\begin{equation}\label{eq:SWME_con}
    \begin{pmatrix}
    1 & & \\ & 1 & \\ & & M^L
\end{pmatrix}\frac{\partial}{\partial t}\begin{pmatrix}
    h \\ h u_m \\ h \alpha
\end{pmatrix} + A^L \frac{\partial}{\partial x}\begin{pmatrix}
    h \\ h u_m \\ h \alpha 
\end{pmatrix} = 0.
\end{equation} 

Since we assume the $N$ spline basis functions have piecewise degree $N$, there is no interior point in the discretization grid and the Legendre basis functions of degree $N$ are in the span of this basis. We can thus represent the Legendre basis functions $\phi^{Leg}_i$ using spline basis functions $\phi^S_i$ via 
\begin{equation}\label{eq:LfromS}
    \phi^{Leg}_i(\zeta) = \sum\limits_{l=1}^{N} T_{il} \phi^S_l(\zeta)
\end{equation}
and a projection of \eqref{eq:LfromS} onto spline basis functions $\phi^S_j$ yields
\begin{equation}
    \int_0^1 \phi^{Leg}_i(\zeta) \phi^S_j(\zeta) d\zeta = \sum\limits_{l=1}^{N} T_{il} \int_0^1 \phi^S_l(\zeta) \phi^S_j(\zeta) d\zeta.
\end{equation}
For the inner products, we use the notation $M^{LS}_{ij} = \int_0^1 \phi^{Leg}_i(\zeta) \phi^S_j(\zeta) d\zeta$ and as before $M^{S}_{ij} = \int_0^1 \phi^S_i(\zeta) \phi^S_j(\zeta) d\zeta$. We then find the transformation matrix $T$ in \eqref{eq:LfromS} from the spline basis to the Legendre basis as $T = M^{LS} (M^S)^{-1}$ so that in vector notation we have
\begin{equation}
    \phi^{Leg} = T \phi^S.
    \label{eq:trafo_basis_SL}
\end{equation}
With the help of the transformation in \eqref{eq:trafo_basis_SL}, we can easily transform the expressions occurring in the system matrix $A^S$ \eqref{eq:A}, i.e., the matrix $M^S$ and the tensors $A^S,B^S$. We obtain after some simple manipulations the transformed versions of the matrices and tensors denoted in the appendix \ref{sec:Matrices} as
\begin{equation}
    M^S = T^{-1} M^L T^{-T}, A_{*j*}^S = (T^{-1} A^L T^T)_{**j} T^{-T}, B_{**j}^S = (T^{-1} B^L T^T)_{**j} T^{-T}.
\end{equation}
The transformation from the spline basis coefficients $s_i$ to the Legendre basis coefficients $\alpha_i$ is similarly done using the velocity profile
\begin{equation}
    u_m + \sum\limits_{i=1}^{N}s_i \phi^S_i(\zeta) = u_m + \sum\limits_{i=1}^{N} s_i \sum\limits_{j=1}^{N} T^{-1}_{ij} \phi^{Leg}_j(\zeta) = u_m + \sum\limits_{j=1}^{N} \underbrace{\sum\limits_{i=1}^{N} s_i T^{-1}_{ij}}_{\alpha_j}\phi^{Leg}_j(\zeta),
\end{equation}
so that we have 
\begin{equation}
    \alpha_1 = T^{-T} \cdot s.
    \label{eq:trafo_coefs_SL}
\end{equation}
Inserting both the transformed basis functions \eqref{eq:trafo_basis_SL} and the transformed coefficients \eqref{eq:trafo_coefs_SL} into the system matrix $A^S$ in \eqref{eq:A}, we obtain from \eqref{eq:SSWME_con} the transformed system 
\begin{equation}
    \begin{pmatrix}
    1 & & \\ & 1 & \\ & & M^{LS} M^L
\end{pmatrix}\frac{\partial}{\partial t}\begin{pmatrix}
    h \\ h u_m \\ h \alpha
\end{pmatrix} + \begin{pmatrix}
    1 & & \\ & 1 & \\ & & M^{LS}
\end{pmatrix} A^L \frac{\partial}{\partial x}\begin{pmatrix}
    h \\ h u_m \\ h \alpha 
\end{pmatrix} = 0.
\end{equation} 

After multiplication with the inverse of the matrix in front of the time derivatives, we obtain
\begin{equation}
    \frac{\partial}{\partial t}\begin{pmatrix}
    h \\ h u_m \\ h \alpha
\end{pmatrix} + \begin{pmatrix}
    1 & & \\ & 1 & \\ & & (M^{L})^{-1}
\end{pmatrix} A^L \frac{\partial}{\partial x}\begin{pmatrix}
    h \\ h u_m \\ h \alpha 
\end{pmatrix} = 0,
\end{equation}
which is exactly the Legendre model or SWME system derived in \cite{KowalskiTorrilhon2019}, where the system matrix is denoted by 
\begin{equation}
   A_{sys}^L= \begin{pmatrix}
    1 & & \\ & 1 & \\ & & (M^{L})^{-1}
\end{pmatrix} A^L. 
\end{equation}
In \cite{KoellermeierRominger2020} it was shown that the Legendre system matrix $A_{sys}^L$ looses hyperbolicity for $N\geq 2$ even for values close to equilibrium. It was shown that the evaluation of the system matrix $A_{sys}^L$ at values $(h,u_m, \tilde{\alpha})=(h, u_m, \alpha_1, 0, \ldots,0)$ always leads to a hyperbolic system for arbitrary $N$. This corresponds to assuming a linear velocity profile and led to a hyperbolic system for the Legendre basis (called HSWME in \cite{KoellermeierRominger2020}), which we denote as
\begin{equation}
    \frac{\partial}{\partial t}\begin{pmatrix}
    h \\ h u_m \\ h \alpha
\end{pmatrix} + A_{sys}^L(h,u_m, \tilde{\alpha}) \frac{\partial}{\partial x}\begin{pmatrix}
    h \\ h u_m \\ h \alpha 
\end{pmatrix} = 0.
\end{equation} 
After the regularization, the transformation back to the spline basis is performed analogously using \eqref{eq:trafo_coefs_SL} and the spline system is then given by 
\begin{equation}\label{eq:retrafoHSSWME}
    \frac{\partial}{\partial t}\begin{pmatrix}
    h \\ h u_m \\ h s
\end{pmatrix} + \begin{pmatrix}
    1 & & \\ & 1 & \\ & & T^{T}
\end{pmatrix} A_{sys}^L(h,u_m, \tilde{s}) \begin{pmatrix}
    1 & & \\ & 1 & \\ & & T^{-T}
\end{pmatrix} \frac{\partial}{\partial x}\begin{pmatrix}
    h \\ h u_m \\ h s 
\end{pmatrix} = 0.
\end{equation} 
Here, $\tilde{s}$ denotes the spline coefficients resulting in zero higher order Legendre coefficients $\tilde{\alpha}$, i.e., $\tilde{\alpha}_i=0$ for $i\geq 2$. This means that $\tilde{s} = T^{T} \tilde{\alpha}$, or more simply 
$\tilde{s}_i= T_{1i} \alpha_1$. Note that since $\tilde{s}$ depends only on the transformation and the first Legendre coefficient $\alpha_1$, it is independent of the higher order coefficients.

The transformed system matrix of the resulting regularized spline model \eqref{eq:retrafoHSSWME} is denoted as $A_{sys}^S(h,hu_m,\tilde{s})$ and given by $diag(1,1,T^T) A_{sys}^L(h,hu_m, \tilde{s}) diag(1,1,T^T)^{-1}$. It is therefore a similarity transformation of the hyperbolic Legendre system matrix $A_{sys}^L(h,hu_m, \tilde{s})$ and has the same eigenvalues, resulting in a hyperbolic spline system matrix. This gives a constructive proof of how a hyperbolic spline system can be derived.

According to the explanations above, the transformations of the system do not have to be performed back and forth. It suffices to first compute from the spline coefficients the projected first Legendre coefficient $\alpha_1$ and 
then compute the linearized spline coefficients $\tilde{s}_i= T_{1i} \alpha_1$. Those linearized spline coefficients $\tilde{s}$ will then be used in the evaluation of the spline system matrix $A_{sys}^S(h,u_m, \tilde{s})$ by
\begin{equation}
    A_{sys}^S = \begin{pmatrix}
    1 & & \\ & 1 & \\ & & M^S
\end{pmatrix}^{-1} A^S,
\end{equation} 
resulting in a hyperbolic system matrix.

We have shown that the regularization around a linear profile leads to hyperbolic systems for all systems based on an $N$ spline basis of degree $N$. This results in the hyperbolic systems HSSWME-L1, HSSWME-Q2, ... . We extend this result using the same methodology without general proof to other spline systems, e.g., HSSWME-L2 and HSSWME-Q3 in the appendix. We conjecture that global hyperbolicity holds for all systems regularized in this way and leave a generalized proof for future work. 

\section{Numerical experiments}
\label{sec:Experiments}
The goal of the numerical experiments is to investigate accuracy and robustness of the proposed SSWME. We therefore simulate a traveling wave over planar bottom topography with periodic boundaries within the domain $ x \in [-1,1]$ until $t_{end}=2$ with similar settings as in \cite{KowalskiTorrilhon2019} and the same finite volume solver therein using a spacial resolution of $n_x=200$ cells. Reference solutions were obtained with a point-wise discretization of the reference system using $N_{\zeta}=60$ points, similar to the ones in \cite{KowalskiTorrilhon2019}.
We denote $\alpha_1$ and $\alpha_2$ for the linear and the quadratic portion of the velocity profile, respectively, see also Section \ref{sec:SplineMomentSystem}. Note that $\alpha_1$ and $\alpha_2$ also correspond to the first two moments of the Legendre systems in \cite{KowalskiTorrilhon2019} and can be computed from the spline coefficients $s_i$ as indicated in Section \ref{sec:constrainedsplines}.

We set the initial water height to
\begin{equation}
\label{eq:initheight}
    h_0(x)=1+\exp (3 \cos (\pi  (x+x_0))-4) \text{ with } x
    _0=1/2
\end{equation} and use a linear initial velocity profile with $u_m=0.25$, $\alpha_1=0.25$ and $\alpha_2=0$ if not stated otherwise. The viscosity is set as $\nu=0.1$ and the slip length as $\lambda=0.1$ 
, which means that the experiments incorporate effects of ground friction and viscosity. The gravitational constant is set to $g=1$.

\subsection{Convergence study: Number of splines}
First, we present a parameter study with increasing number of splines $N$ or equivalently increasing resolution of the vertical grid, to validate the SSWME models.
\begin{figure}
    \label{fig:spline-number}
    \centering
    \includegraphics[width=0.48\textwidth]{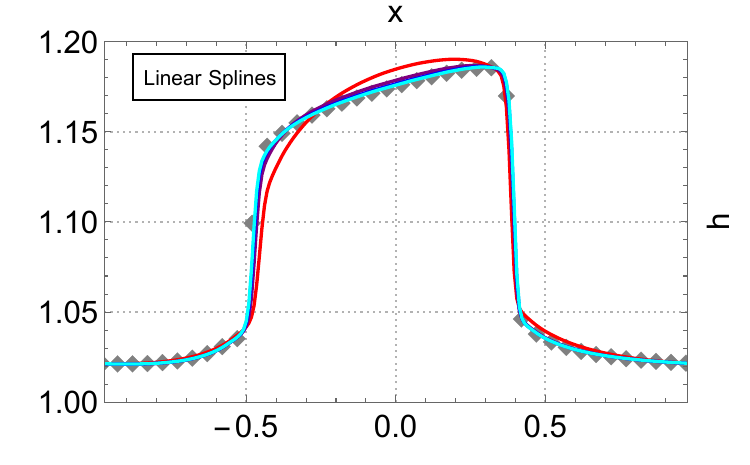}
    \includegraphics[width=0.48\textwidth]{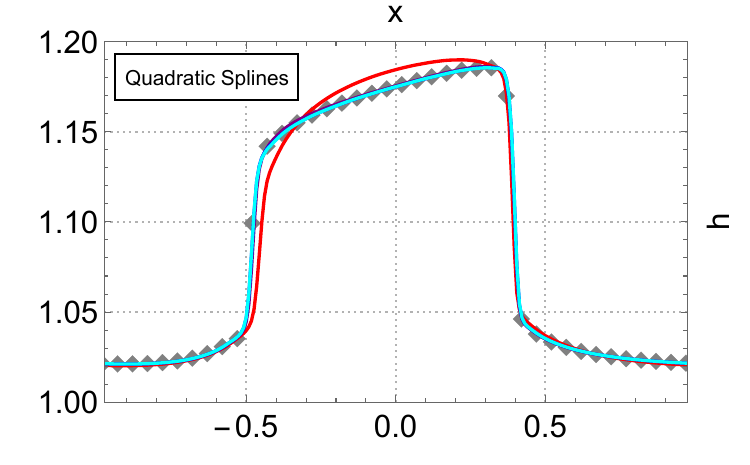}\\
    \includegraphics[width=0.48\textwidth]{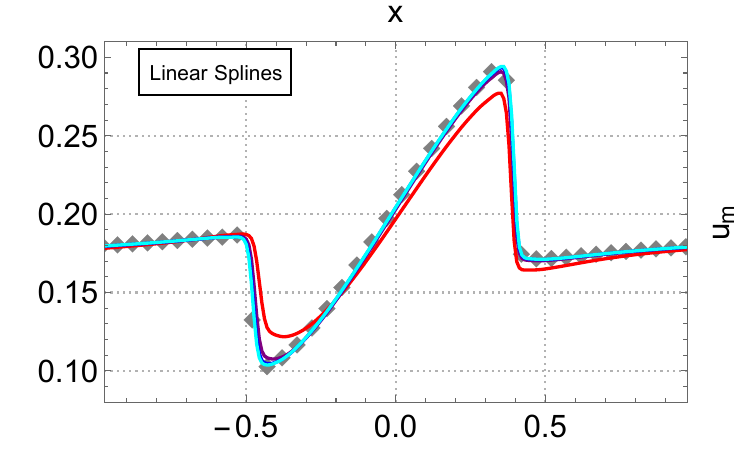}
    \includegraphics[width=0.48\textwidth]{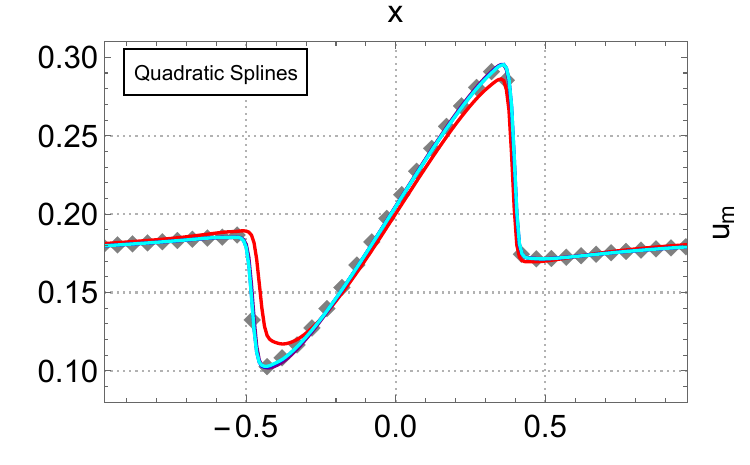}\\
    \includegraphics[width=0.48\textwidth]{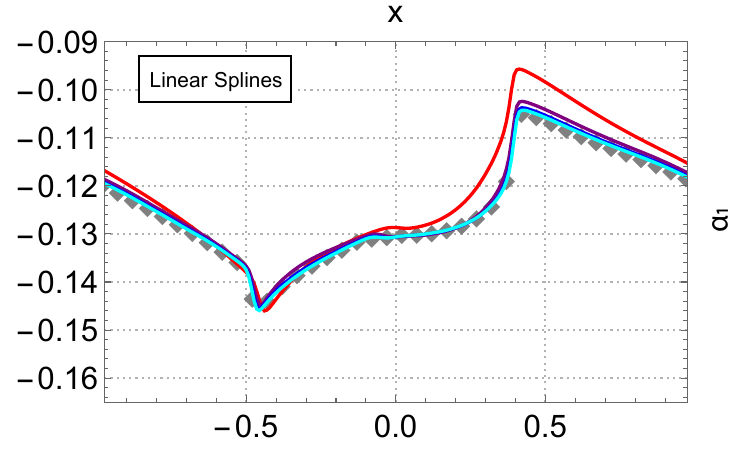}
    \includegraphics[width=0.48\textwidth]{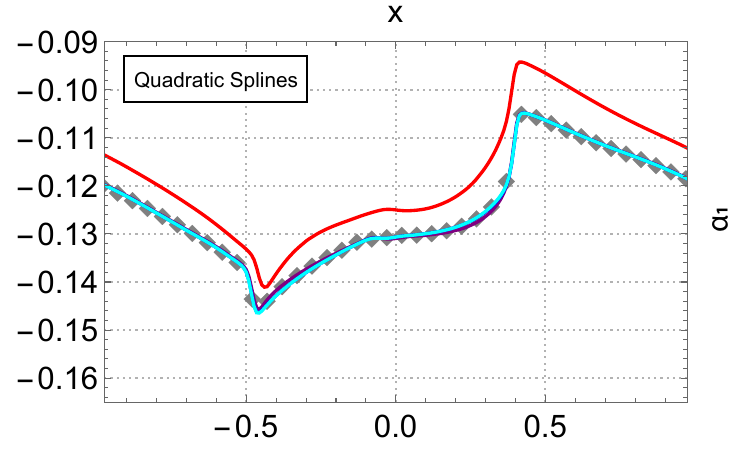}\\
    \includegraphics[width=0.48\textwidth]{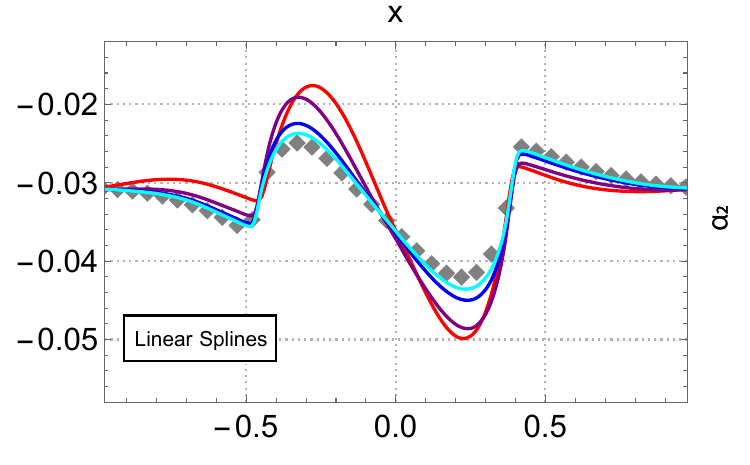}
    \includegraphics[width=0.48\textwidth]{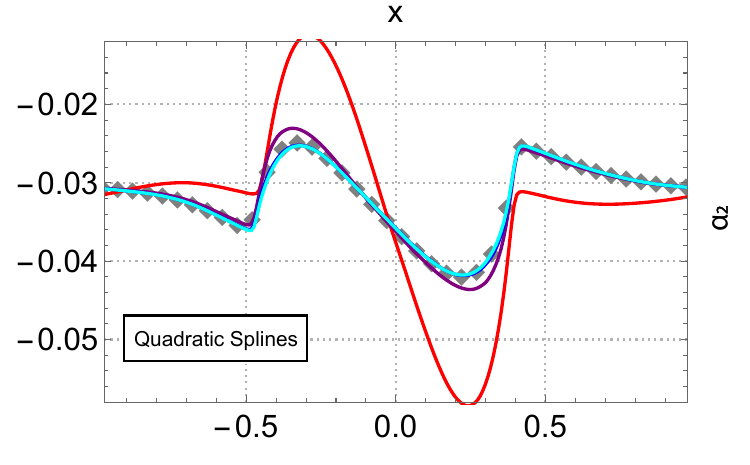}\\
    \hspace{2ex}\includegraphics[width=.3\textwidth]{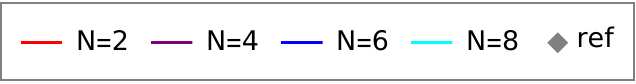}
    \caption{Results of the smooth wave experiment with  N=2 (red) N=4 (purple) N=6 (blue) and N=8 (cyan) spline moments for the variables height $h$ (top row), mean velocity $u_m$ (second row), linear portion $\alpha_1$ (third row), quadratic portion $\alpha_2$ (last row). Left column: linear systems SSWME-LN. Right column: quadratic systems SSWME-QN. All models converge with increasing number of moments $N$.
    }
    \label{fig:numberOfSplines}
\end{figure}

\begin{figure}
    \centering
    \begin{overpic}[width=0.48\linewidth]{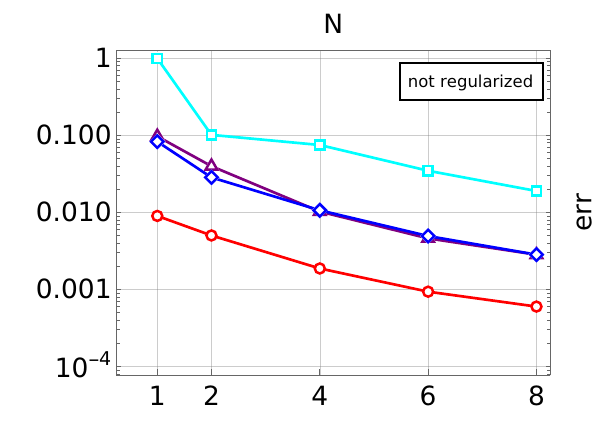}
    \put(25,29){$\ured{h}$}
    \put(25,42){$\ublue{\alpha_1}$}
    \put(25,50){$\upurple{u_m}$}
    \put(25,55){$\ucyan{\alpha_2}$}
    \end{overpic}
    \begin{overpic}[width=0.48\linewidth]{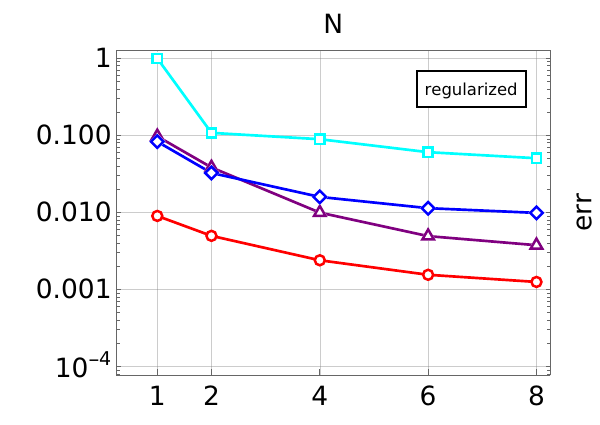}
    \put(25,29){$\ured{h}$}
    \put(25,42){$\ublue{\alpha_1}$}
    \put(25,50){$\upurple{u_m}$}
    \put(25,55){$\ucyan{\alpha_2}$}
    \end{overpic}\\
    \begin{overpic}[width=0.48\linewidth]{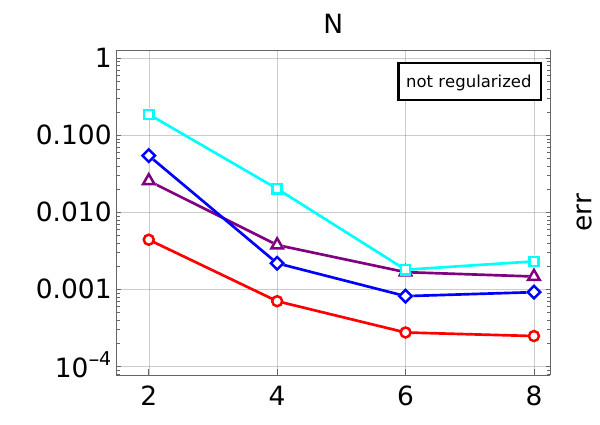}
    \put(25,25){$\ured{h}$}
    \put(25,47){$\ublue{\alpha_1}$}
    \put(25,35){$\upurple{u_m}$}
    \put(25,55){$\ucyan{\alpha_2}$}
    \end{overpic}
    \begin{overpic}[width=0.48\linewidth]{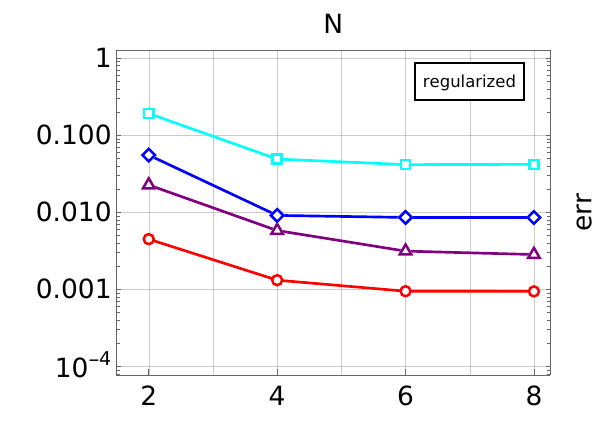}
    \put(25,25){$\ured{h}$}
    \put(25,47){$\ublue{\alpha_1}$}
    \put(25,35){$\upurple{u_m}$}
    \put(25,55){$\ucyan{\alpha_2}$}
    \end{overpic}
    \\
    \begin{overpic}[width=0.48\linewidth]{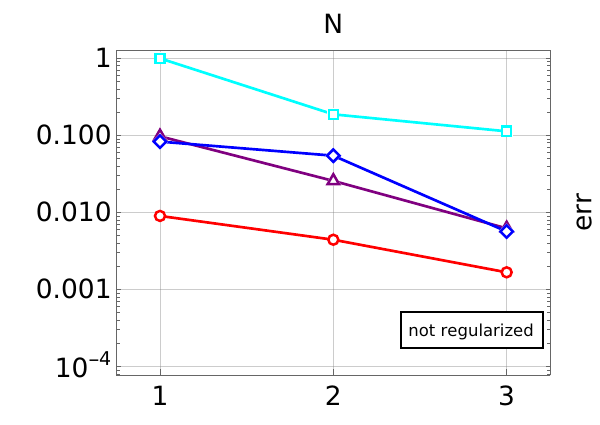}
       \put(25,29){$\ured{h}$}
    \put(25,42){$\ublue{\alpha_1}$}
    \put(25,50){$\upurple{u_m}$}
    \put(25,55){$\ucyan{\alpha_2}$}
        \end{overpic}
        \begin{overpic}[width=0.48\linewidth]{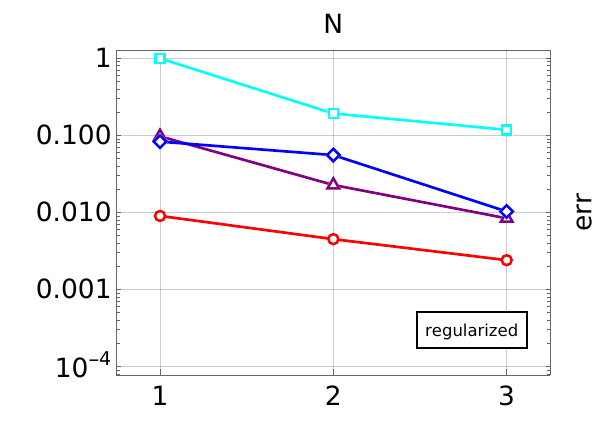}
           \put(25,29){$\ured{h}$}
    \put(25,42){$\ublue{\alpha_1}$}
    \put(25,50){$\upurple{u_m}$}
    \put(25,56){$\ucyan{\alpha_2}$}
        \end{overpic}
    \caption{Relative $L_1$ error $err:=\frac{||v-v_{ref}||_1}{||v_{ref}||_1}$ with $v$ either water height $h$ (red), mean velocity $u_m$ (purple), linear coefficient $\alpha_1$ (blue), or quadratic coefficient $\alpha_2$ (cyan) for different numbers of splines $N$; linear spline models SSWME-LN (top row), quadratic spline models SSME-QN (middle row); standard Legendre models SWMEN (bottom row); non-hyperbolic models (left column), regularized hyperbolic models (right column). The error decreases the more splines are used. The quadratic spline systems SSWME-QN show lower errors than the corresponding linear systems SSWME-LN.}
    \label{fig:spline-number-error}
\end{figure}

From Figure \ref{fig:numberOfSplines} we observe that the model accuracy increases with increasing number of moments $N$. This is true both for the linear spline models SSWME-LN as well as for the quadratic spline models SSWME-QN. The behavior of $\alpha_1$ and $\alpha_2$ seems to suggest that the quadratic spline models SSWME-QN converge faster than the linear spline models SSWME-LN. 
However, we note that the quadratic system with only two splines SSWME-Q2 shows a larger deviation in $\alpha_1$ and $\alpha_2$ than the corresponding linear system SSWME-L2. This is due to the flexibility and regularity restriction of the two quadratic splines which are $C^2$, as compared to the $C^1$ linear splines with the additional grid point at $\zeta=\frac{1}{2}$.

As shown in Figure \ref{fig:spline-number-error} (left), the approximations obtained by linear models SSWME-LN and the quadratic models \mbox{SSWME-QN} for increasing number of splines $N$ from from 1 to 8 both converge to the reference solution as the relative $L_1$ errors in water height $h$, mean velocity $u_m$, $\alpha_1$, and $\alpha_2$ decrease the more splines are used.
For the quadratic spline models SSWME-QN the error saturates at $N=6$. This can be attributed to remaining numerical errors in the treatment of the friction terms or numerical diffusion in the solvers provided by \cite{KowalskiTorrilhon2019}. However, as expected the quadratic spline systems achieve a higher accuracy than the linear spline systems with the same number of splines.

\subsection{Effect of the hyperbolic regularization}
\begin{figure}
    \centering
    \includegraphics[width=0.4\textwidth]{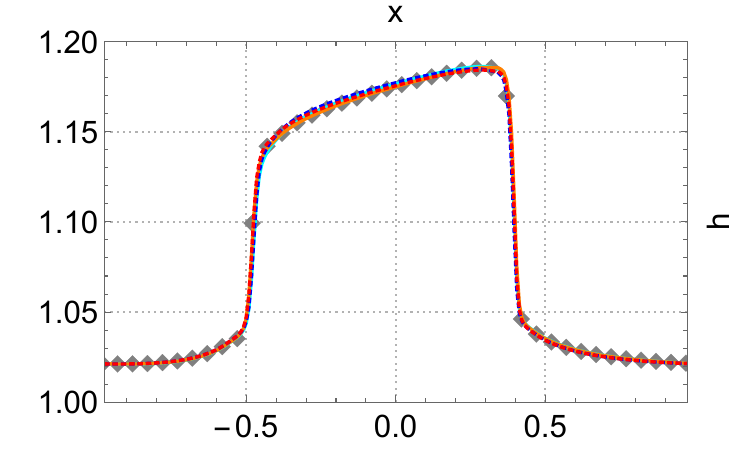}
    \includegraphics[width=0.4\textwidth]{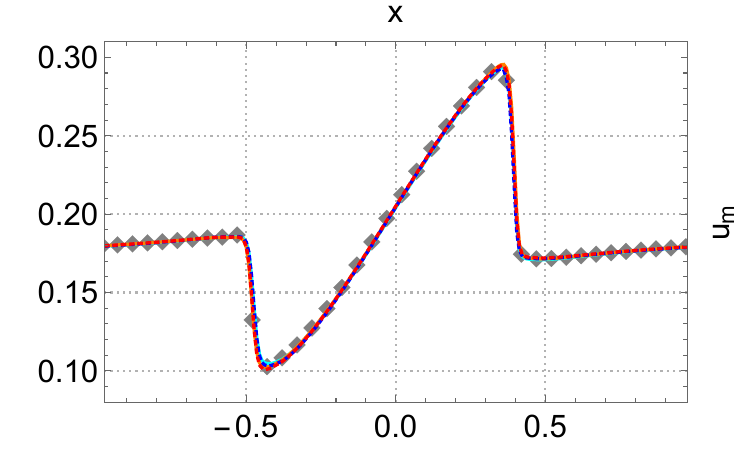}\\
    \includegraphics[width=0.4\textwidth]{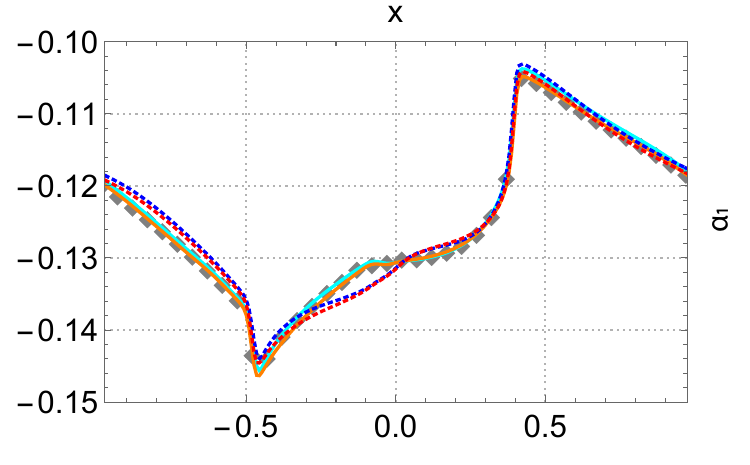}
    \includegraphics[width=0.4\textwidth]{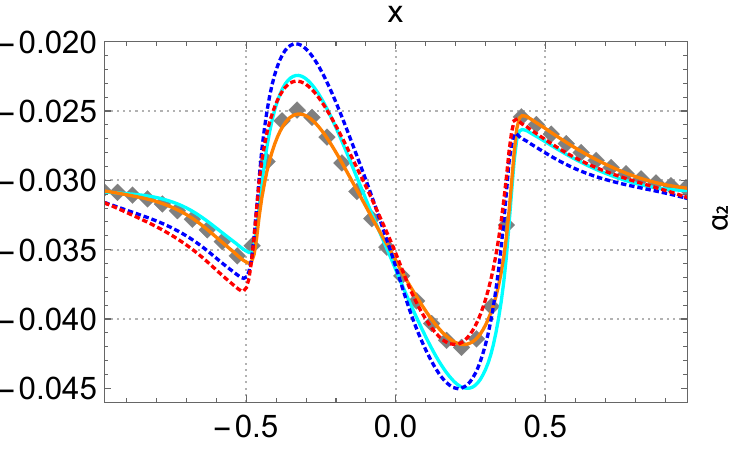}\\
    \includegraphics[width=.5\textwidth]{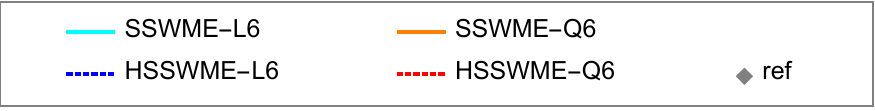}
    \caption{Results of the smooth wave experiment obtained with the models SSWME-L6 (cyan) and SSWME-Q6 (orange) and their regularized versions HSSWME-L6 (blue) and HSSWME-Q6 (red) for the quantities water height $h$ and mean velocity $u_m$ (top row) and linear and quadratic profile portions $\alpha_1$ and $\alpha_2$ (second row). The velocity profile is less accurately modeled by the regularized equations for some regions.}
    \label{fig:reg-compare}
\end{figure}

As shown in Figure \ref{fig:spline-number-error} (right), the observed convergence with respect to the spline number is mostly the same in the hierarchy of the regularized systems HSSWME-LN and HSSWME-QN, too. As a consequence of an additional model error introduced by the regularization, the error saturation occurs slightly earlier. 

Figure \ref{fig:reg-compare} shows how the numerical solution of the regularized systems behaves compared to the non-regularized systems. In the second row the regularized systems show some discrepancies in the approximation of the higher moments of the velocity profile $\alpha_1$ and $\alpha_2$. However, water height $h$ and mean velocity $u_m$ show visually no differences.

\subsection{Comparison with classical shallow water moment equations}
\begin{figure}
    \centering
    \includegraphics[width=0.4\textwidth]{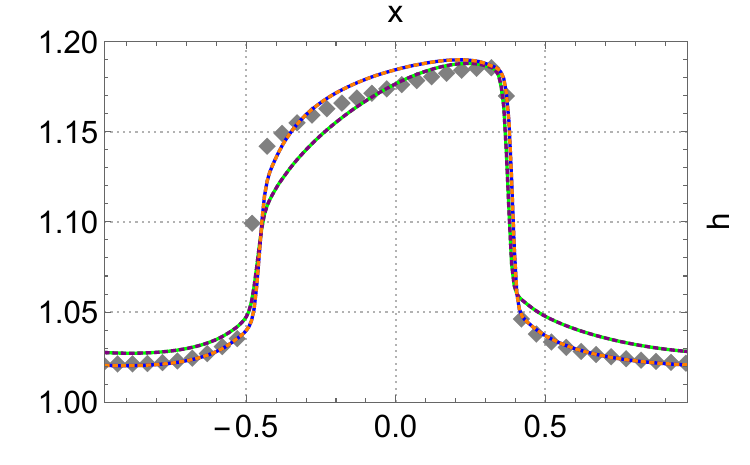}
    \includegraphics[width=0.4\textwidth]{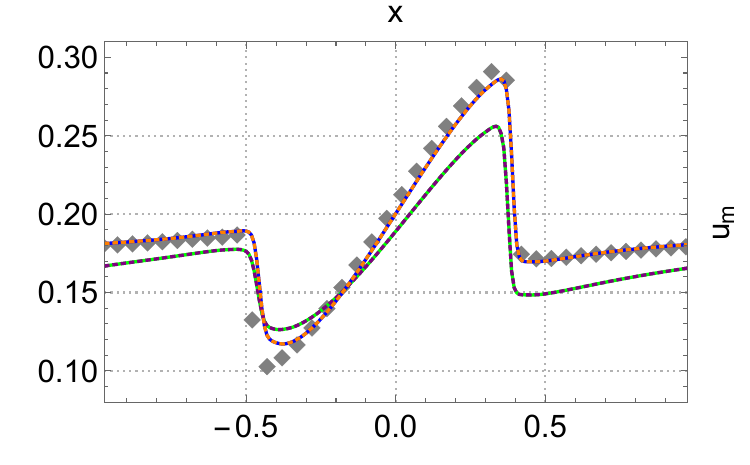}\\
    \includegraphics[width=0.4\textwidth]{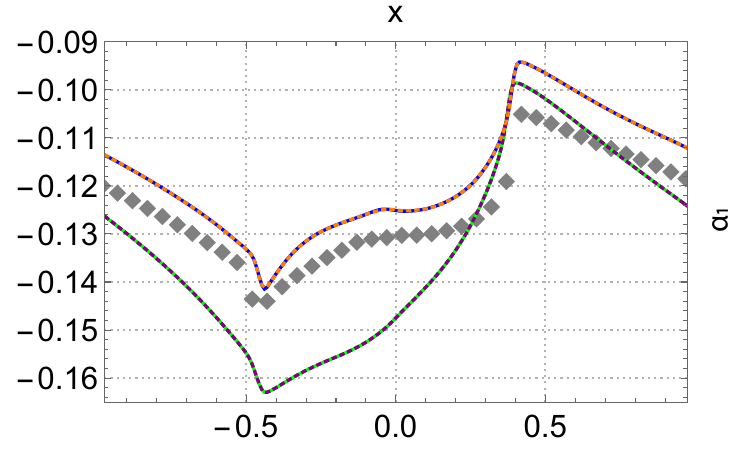}
    \includegraphics[width=0.4\textwidth]{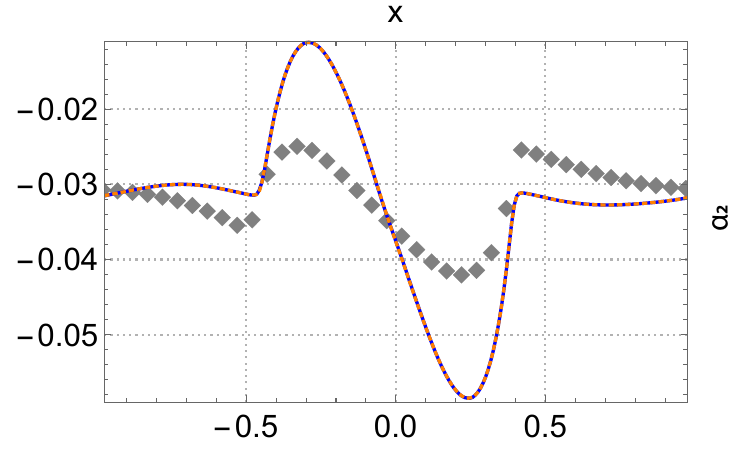}\\
    \includegraphics[width=.45\textwidth]{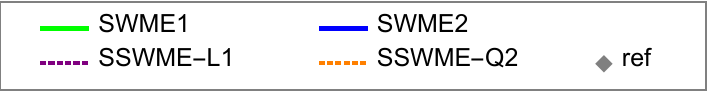}
    \caption{Results of the smooth wave experiment using the spline systems SSWME-L1 (purple, dashed) and SSWME-Q2 (orange, dashed) and the Legendre models from the shallow water moments hierarchy SWME1 (green) and SWME2 (blue). Depicted are the height $h$ (top left), velocity $u_m$ (top right) and the linear and quadratic profile portions $\alpha_1$ (bottom left) and $\alpha_2$ (bottom right). The curves of equivalent systems align.}
    \label{fig:sme-equal}
\end{figure}

\begin{figure}
    \centering
    \includegraphics[width=0.4\textwidth]{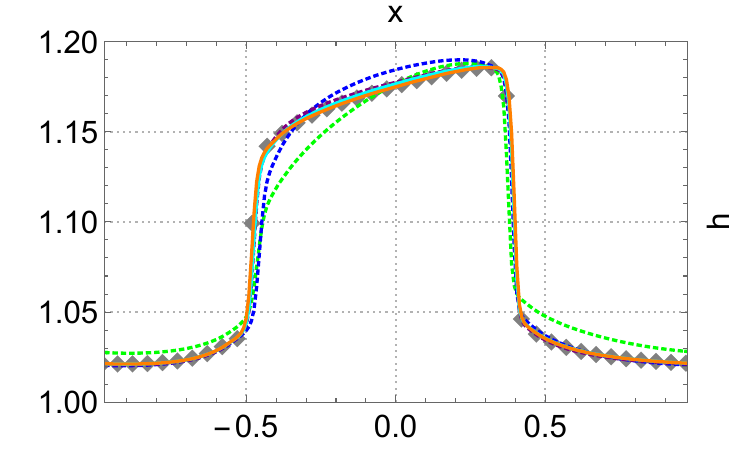}
    \includegraphics[width=0.4\textwidth]{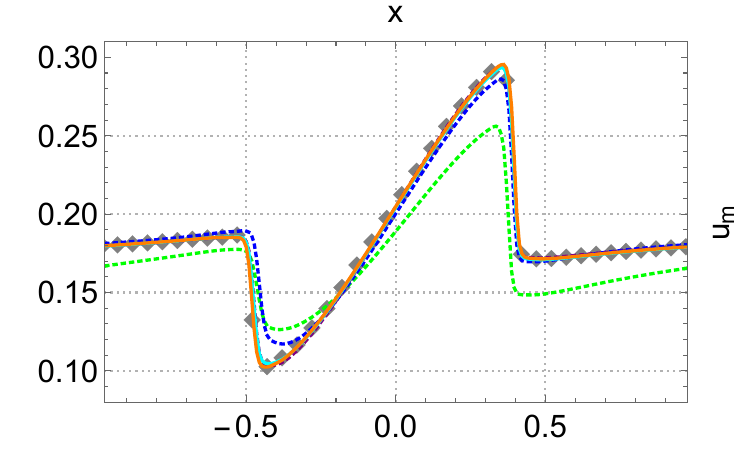}\\
    \includegraphics[width=0.4\textwidth]{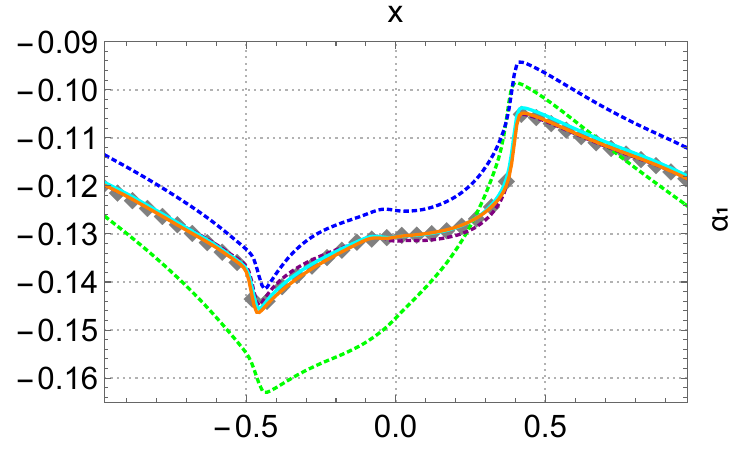}
    \includegraphics[width=0.4\textwidth]{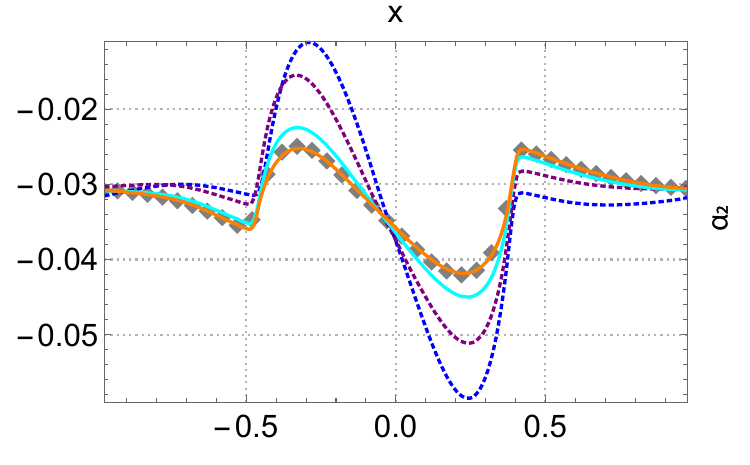}\\
    \includegraphics[width=.45\textwidth]{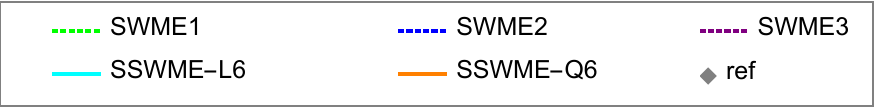}
    \caption{Results of the smooth wave experiment obtained with spline models SSWME-L6 (cyan) and SSWME-Q6 (orange) and Legendre models SWME1 (green), SWME2 (blue) and SWME3 (red). The plots show water height $h$ and mean velocity $u_m$ (top row) and linear and quadratic profile portions $\alpha_1$ and $\alpha_2$ (bottom row). Especially the higher order velocity moments $\alpha_1$ and $\alpha_2$ are best approximated by the spline systems, with SSWME-Q6 leading to the smallest deviation from the reference solution.}
    \label{fig:sme-compare}
\end{figure}

The results obtained with the SSWME are consistent with the known results for the SWME from \cite{KowalskiTorrilhon2019} in their non-regularized and regularized form. 

Figure \ref{fig:sme-equal} presents the numerical proof of the equality explained earlier between the systems SWME1 and SSWME-L2, and the systems SWME2 and SSWME-Q2, respectively. 

Figure \ref{fig:sme-compare}, on the other hand, shows that it is possible to achieve more accurate results with models from the SSWME hierarchy such as SSWME-L6 and SSWME-Q6 in comparison with the known systems SWME1 to SWME3, which were the most accurate systems implemented in the software accompanying \cite{KowalskiTorrilhon2019}.

\subsection{Fast wave with steep initial profile}
\begin{figure}
    \centering
    \includegraphics[width=0.46\linewidth]{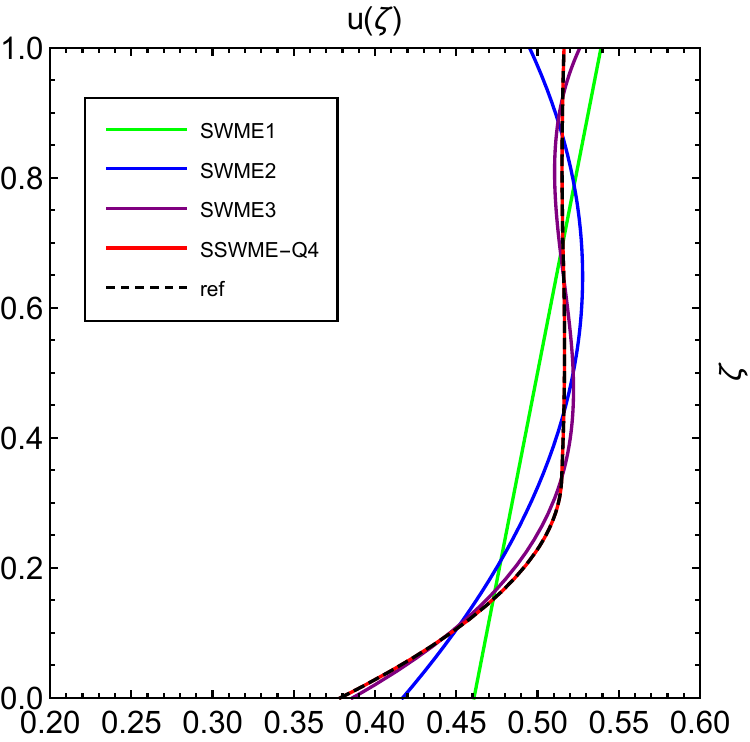}
    \caption{Approximations of the initial velocity profile function \eqref{eq:IC} in the fast wave experiment with linear (green), quadratic (blue), and cubic (purple) Legendre polynomials compared with $N=4$ quadratic splines Q4 (red). The piecewise defined initial profile is only approximated well by the spline basis.}
    \label{fig:steep-init-profile}
\end{figure}

In a different setting we employ a fast and steep initial velocity profile $u(x,0,\zeta)$ realized piecewise and continuously by
\begin{equation}\label{eq:IC}
    u(x,0,\zeta):=\begin{cases}
\begin{array}{cc}
 \frac{15 \left(1569 \xi ^2-2568 \xi +1700\right)}{322048} & \frac{2}{3}\leq \xi \leq 1, \\
 -\frac{15 \left(1887 \xi ^2-2040 \xi -164\right)}{322048} & \frac{1}{3}\leq \xi <\frac{2}{3}, \\
 -\frac{15 \left(50298 \xi ^2-34314 \xi +5215\right)}{322048} & 0\leq \xi <\frac{1}{3}, \\
\end{array}
\end{cases}
\end{equation}
as depicted in Figure \ref{fig:steep-init-profile}. This means a faster background flow with average velocity $u_m=0.5$. In addition, we consider almost kinematic conditions with viscosity $\nu=0.0005$ and slip length $\lambda=0.1$. The initial height function remains as defined in \eqref{eq:initheight}.

The SSWME systems are better suited than the SWME for flexibly modeling profile shapes during the simulation. Figure \ref{fig:steep-init-profile} shows the approximation of the initial condition profile \eqref{eq:IC} where the near-bottom region resembles the steep increase of a logarithmic profile, typically encountered in boundary layers. The SSWME's Q4 basis approximates this profile exactly, as it employs a piecewise quadratic basis that includes \eqref{eq:IC} in the approximation space. On the other hand, oscillations and deviations occur in the approximations with the Legendre systems. Note that these oscillations do not vanish with increasing Legendre degree, as the piecewise continuous profile is not in the space of global polynomials on $[0,1]$. 

\begin{figure}
    \centering
    \includegraphics[width=0.48\linewidth]{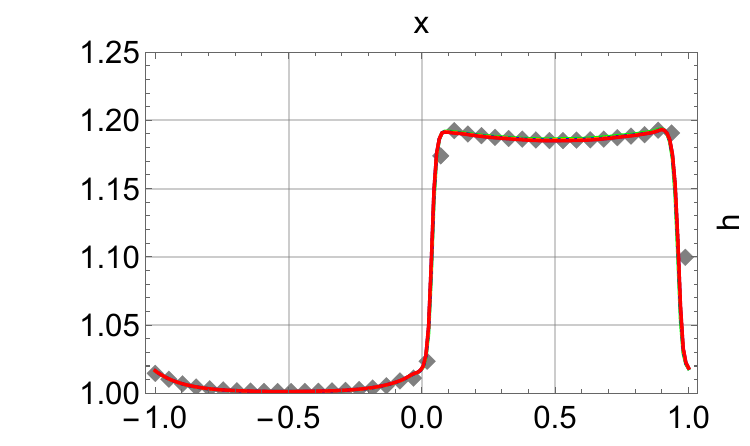}
    \includegraphics[width=0.48\textwidth]{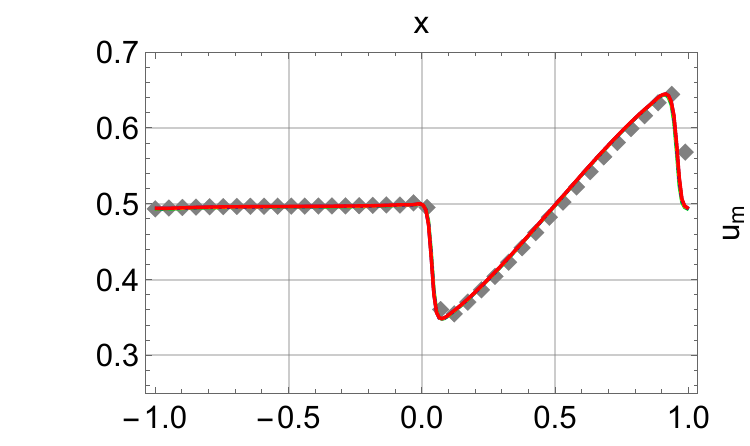} \\
    \includegraphics[width=0.48\textwidth]{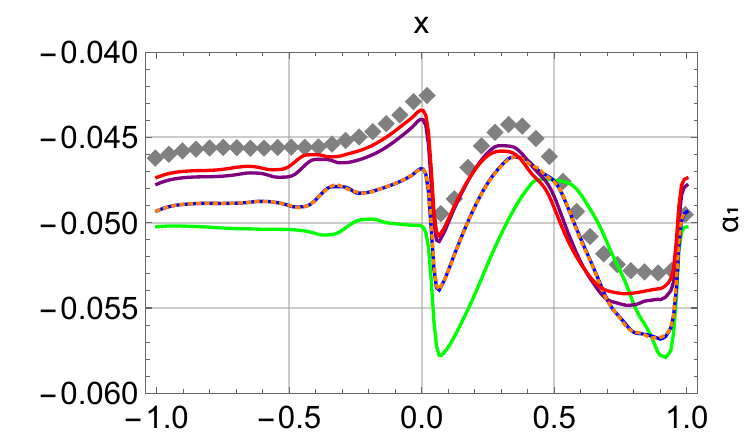}
    \includegraphics[width=0.48\textwidth]{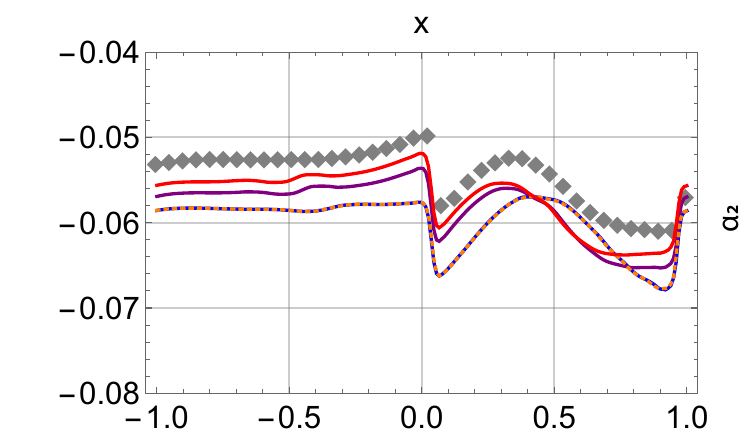}\\
    \includegraphics[width=0.5\linewidth]{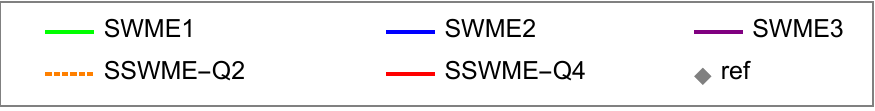}
    \caption{Results of the fast wave experiment using the Legendre models SWME1 (green), SWME2 (blue) and SWME3 (purple), and the quadratic spline models SSWME-Q2 (orange) and SSWME-Q4 (red). The plots show the water height $h$ and mean velocity $u_m$ (first row) as well as profile portions $\alpha_1$ and $\alpha_2$ (second row). Due to the bottom velocities being too large and the resulting difference in propagation speeds, all systems slightly deviate from the velocity of the wave as the lines indicating height $h$ and velocity $u_m$ lie to the left of the reference solution. However, SSWME-Q4 yields the best approximation.}
    \label{fig:steep-h-um}
\end{figure}

\begin{figure}
    \centering
    \includegraphics[width=0.46\linewidth]{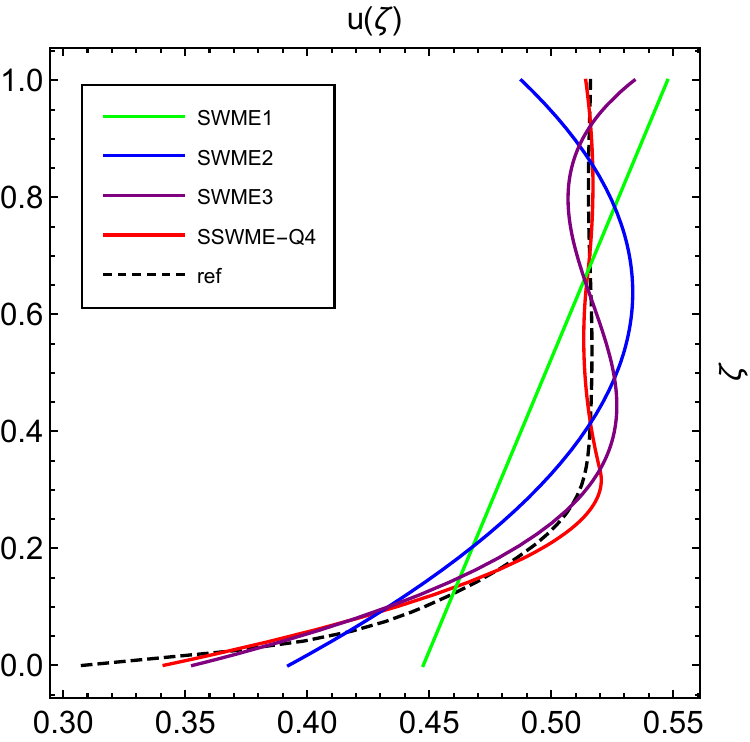}
    \includegraphics[width=0.46\linewidth]{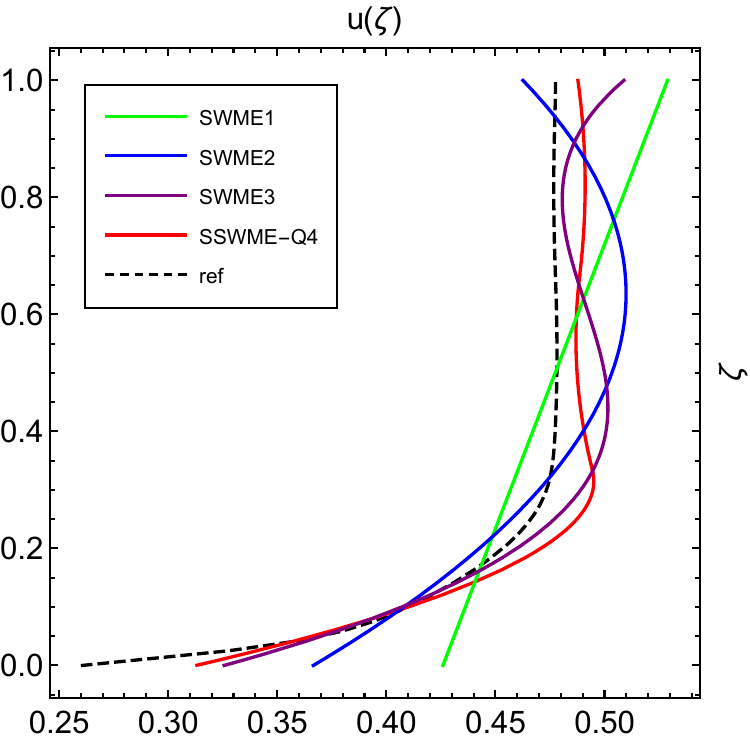}
    \caption{Predicted velocity profiles in the fast wave experiment using the Legendre models SWME1 (green), SWME2 (blue), and SWME3 (purple), and the spline model SSWME-Q4 (red). Left: $x=-0.105$. Right: $x=0.025$. The profile has changed due to bottom friction. The Legendre models are less accurate than SSWME-Q4. All models overestimate the bottom velocity, leading to different propagation speeds. 
    }
    \label{fig:steep-profiles}
\end{figure}

Figure \ref{fig:steep-h-um} shows the transient simulation results obtained using the steep profile \eqref{eq:IC}. The SSWME-Q4 system performs better than the SSWME-Q2 system and the SWME systems which is best visible in the approximation of the quantities $\alpha_1$ and $\alpha_2$. Note that all systems slightly deviate from the velocity of
the wave as the lines indicating height $h$ and velocity $u_m$ lie to the left of the reference solution. This could be attributed to slightly wrong bottom velocities as seen in Figure \ref{fig:steep-profiles} for all models, in turn leading to larger bottom friction and smaller $\alpha_1$, $\alpha_2$ values as well as propagation speeds. From this it becomes clear that friction coefficients, which are typically calibrated for depth-averaged models only, need to be remodeled. However, this is beyond the scope of this paper. 

In Figure \ref{fig:steep-profiles}, the velocity profiles at $t=2$ are shown at two points in the domain. It is clear that the profile shape changes during the simulation due to friction at the ground. We observe that the SSWME-Q4 model outperforms the Legendre models with respect to the accuracy of the velocity profile.

In both our test cases the SSWME remained inside the hyperbolic parameter range and we did not encounter problems with instability. Wherever instabilities occur, HSSWME might overcome them as shown in \cite{KoellermeierRominger2020} for the SME. We leave this analysis for future work.

\section{Conclusion and future work}
This paper introduces flexible Spline Shallow Water Moment Equations (SSWME) as reduced models for free-surface flows. The models are derived after a construction of appropriate constrained spline basis functions using a Galerkin projection. The closed form of the reduced models allows for in-depth analysis. We show that hyperbolicity loss can be mitigated using a regularization approach and give a constructive hyperbolicity proof for a hierarchy of hyperbolic models of arbitrary order. 
The numerical experiments show that the SSWME are a suitable means for dimensional reduction of free surface flow preserving information of the velocity profile. As was shown, an additional benefit of the spline basis is the ability to approximate velocity profiles with less regularity or steep gradients.


A possibility for future work is the use of irregular spline grids or adaptive bases such as wavelets. Furthermore, other hyperbolic regularizations and a generalized proof could lead to the development of new model hierarchies. The inclusion and calibration of more physical parameters together with more numerical simulations would increase the applicability of the derived models. 

\section*{Acknowledgments}
This publication is part of the project \textit{HiWAVE} with file number VI.Vidi.233.066 of the \textit{ENW Vidi} research programme, funded by the \textit{Dutch Research Council (NWO)} under the grant \url{https://doi.org/10.61686/CBVAB59929}.

\appendix
\begin{section}{Appendix}
\subsection{Vectors, matrices and tensors in the general moment equations}
\label{sec:Matrices}
The vectors, matrices, and tensors in the general form of the moment equations \eqref{eq:momenteqn} are
    \begin{equation}
        M_{ij}=\int_0^1 \phi_i(\zeta) \phi_j(\zeta) \text{ d}\zeta,
        \qquad
        V^0_i=\phi_i(0),
        \qquad
        C_{ij}=\int_0^1 \phi_i'(\zeta) \phi_j'(\zeta) \text{ d}\zeta,
    \end{equation}
    \begin{equation}
        A_{ijk}=\int_0^1 \phi_i(\zeta)\phi_j(\zeta)\phi_k(\zeta)\text{ d}\zeta,
    \qquad
        B_{ijk}=\int_0^1 \phi_i'(\zeta)
        \int_0^{\zeta}\phi_j(\hat\zeta)\text{d}\hat\zeta\,\phi_k(\zeta)\text{ d}\zeta.
    \end{equation}

\subsection{SSWME system matrices and characteristic polynomials}
\label{sec:systems}
We give an overview of the SSWME models included in this paper for linear and quadratic splines and $N=1, 2, 3$ by providing the system matrices $A_{sys}$, the characteristic polynomials and the hyperbolic regularizations.

\subsubsection{SSWME-L1}
    The SSWME-L1 system matrix is given by
    \begin{equation}
    A_{sys}^{L1} = \begin{pmatrix}
 0 & 1 & 0 \\
 gh-u_m^2-\frac{4 s_1^2}{3} & 2 u_m & \frac{8 s_1}{3} \\
 -2 s_1 u_m & 2 s_1 & u_m \\
\end{pmatrix}
\end{equation}
with shifted characteristic polynomial ($\lambda=u_m + c \sqrt{gh}, s_1=\overline{s_1}\sqrt{gh}$)
\begin{equation}
       c^3-c(4 \bar s_1^2+1)  =0
\end{equation}

The eigenvalues of the system matrix are 
\begin{equation}
    \lambda_{1,2,3} \in \left\{u_m,u_m\pm\sqrt{gh+4 s_1^2}\right\}.
\end{equation}

\subsubsection{(H)SSWME-L2}
The SSWME-L2 system matrix \eqref{eq:ASysSSWME-L2} yields the shifted characteristic polynomial \eqref{eq:char_poly}.
Not all roots of this polynomial are real.
Restricting the coefficients $(s_1,s_2)=(\frac{1}{4},\frac{1}{4})\alpha_1$ yields the system matrix  \eqref{eq:ASysHSSWME-L2} of the regularized HSSWME-L2 system
and the requiring roots of the shifted polynomial \eqref{eq:charpolyH}.
The real eigenvalues of the hyperbolic HSSWME-L2 are then given by 
\begin{equation}
    \lambda_{1,2,3,4} \in \{u_m \pm \frac{\sqrt{3} }{4} \alpha_1 ,u_m \pm \sqrt{ \alpha_1 ^2+gh}\}.
\end{equation}

\subsubsection{(H)SSWME-L3}
The SSWME-L3 system matrix is given by $A_{sys}^{L3} =$
    \begin{multline}
    \left(
        \begin{array}{c c c : }
 0 & 1 & 0 \\
 gh-u_m^2-4 s_1^2+s_2 s_1+s_3 s_1-3 s_2^2-4 s_3^2+s_2 s_3 & 2 u_m & 8 s_1-s_2-s_3 \\
 -2 s_1 u_m-\frac{77 s_1^2}{30}-\frac{26 s_2 s_1}{15}-\frac{s_3 s_1}{3}+\frac{s_2^2}{2}+\frac{11 s_3^2}{15}+\frac{s_2 s_3}{15} & 2 s_1 & u_m+\frac{49 s_1}{12}-\frac{s_2}{24}+\frac{49 s_3}{120} \\
 -2 s_2 u_m-\frac{9 s_1^2}{5}+\frac{6 s_2 s_1}{5}+\frac{9 s_3^2}{5}-\frac{6 s_2 s_3}{5} & 2 s_2 & \frac{9 s_1}{2}-\frac{9 s_2}{8}-\frac{9 s_3}{40} \\
 -2 s_3 u_m-\frac{11 s_1^2}{15}-\frac{s_2 s_1}{15}+\frac{s_3 s_1}{3}-\frac{s_2^2}{2}+\frac{77 s_3^2}{30}+\frac{26 s_2 s_3}{15} & 2 s_3 & \frac{7 s_1}{6}+\frac{s_2}{24}-\frac{31 s_3}{120} \\
\end{array}
\right. \\
\left.
\begin{array}{: c c}
 0 & 0 \\
 -s_1+6 s_2-s_3 & -s_1-s_2+8 s_3 \\
 \frac{169 s_1}{120}-\frac{7 s_2}{4}+\frac{61 s_3}{120} & \frac{31 s_1}{120}-\frac{s_2}{24}-\frac{7 s_3}{6} \\
 u_m+\frac{21 s_1}{40}-\frac{21 s_3}{40} & \frac{9 s_1}{40}+\frac{9 s_2}{8}-\frac{9 s_3}{2} \\
 -\frac{61 s_1}{120}+\frac{7 s_2}{4}-\frac{169 s_3}{120} & u_m-\frac{49 s_1}{120}+\frac{s_2}{24}-\frac{49 s_3}{12} \\
\end{array} \right)
    \end{multline}
    The shifted characteristic polynomial following the notation from \eqref{eq:evform} is
    \begin{multline}
       160c^5-672c^4T_4-2c^3T_3+3c^2T_2-6cT_1+9T_0\\
    \end{multline}
    where \begin{align}
        T_4&=\overline{s_1}-\overline{s_3},\\
        T_3&=1461 \overline{s_1}{}^2-1029 \overline{s_2} \overline{s_1}+1617 \overline{s_3} \overline{s_1}+1035 \overline{s_2}{}^2+1461 \overline{s_3}{}^2-1029 \overline{s_2} \overline{s_3}+80\\
        T_2&=672 \overline{s_1}{}^3-3093 \overline{s_2} \overline{s_1}{}^2-4125 \overline{s_3} \overline{s_1}{}^2+3987 \overline{s_2}{}^2 \overline{s_1}+4125 \overline{s_3}{}^2 \overline{s_1}+224 \overline{s_1}-672 \overline{s_3}{}^3\\&\nonumber+3093 \overline{s_2} \overline{s_3}{}^2-3987 \overline{s_2}{}^2 \overline{s_3}-224 \overline{s_3}\\
        T_1&=252 \overline{s_1}{}^4+45 \overline{s_2} \overline{s_1}{}^3-2763 \overline{s_3} \overline{s_1}{}^3-2628 \overline{s_2}{}^2 \overline{s_1}{}^2-7191 \overline{s_3}{}^2 \overline{s_1}{}^2+12897 \overline{s_2} \overline{s_3} \overline{s_1}{}^2\\&\nonumber-167 \overline{s_1}{}^2+522 \overline{s_2}{}^3 \overline{s_1}-2763 \overline{s_3}{}^3 \overline{s_1}+12897 \overline{s_2} \overline{s_3}{}^2 \overline{s_1}+263 \overline{s_2} \overline{s_1}-12762 \overline{s_2}{}^2 \overline{s_3} \overline{s_1}\\&\nonumber-619 \overline{s_3} \overline{s_1}+135 \overline{s_2}{}^4+252 \overline{s_3}{}^4+45 \overline{s_2} \overline{s_3}{}^3-105 \overline{s_2}{}^2-2628 \overline{s_2}{}^2 \overline{s_3}{}^2-167 \overline{s_3}{}^2\\&\nonumber+522 \overline{s_2}{}^3 \overline{s_3}+263 \overline{s_2} \overline{s_3}\\
        T_0&=180 \overline{s_2} \overline{s_1}{}^4-1548 \overline{s_3} \overline{s_1}{}^4-81 \overline{s_2}{}^2 \overline{s_1}{}^3-1737 \overline{s_3}{}^2 \overline{s_1}{}^3+1638 \overline{s_2} \overline{s_3} \overline{s_1}{}^3\\&\nonumber+32 \overline{s_1}{}^3-774 \overline{s_2}{}^3 \overline{s_1}{}^2+1737 \overline{s_3}{}^3 \overline{s_1}{}^2+7 \overline{s_2} \overline{s_1}{}^2+8721 \overline{s_2}{}^2 \overline{s_3} \overline{s_1}{}^2\\&\nonumber+95 \overline{s_3} \overline{s_1}{}^2+351 \overline{s_2}{}^4 \overline{s_1}+1548 \overline{s_3}{}^4 \overline{s_1}-1638 \overline{s_2} \overline{s_3}{}^3 \overline{s_1}-177 \overline{s_2}{}^2 \overline{s_1}-8721 \overline{s_2}{}^2 \overline{s_3}{}^2 \overline{s_1}\\&\nonumber-95 \overline{s_3}{}^2 \overline{s_1}-180 \overline{s_2} \overline{s_3}{}^4+81 \overline{s_2}{}^2 \overline{s_3}{}^3-32 \overline{s_3}{}^3+774 \overline{s_2}{}^3 \overline{s_3}{}^2-7 \overline{s_2} \overline{s_3}{}^2-351 \overline{s_2}{}^4 \overline{s_3}+177 \overline{s_2}{}^2 \overline{s_3}
    \end{align}

    A hyperbolic regularization around a linear profile using the restriction $s_1=\frac{\alpha_1}{6}$,
    $s_2=\frac{5}{18}\alpha_1$ and again $s_3=\frac{\alpha_1}{6}$ yields a new regularized HSSWME-L3 system
    matrix
\begin{equation}
 A_{sys,H}^{L3} = \begin{pmatrix}
 0 & 1 & 0 & 0 & 0 \\
 -\frac{\alpha_1 ^2}{3}+gh-u_m^2 & 2 u_m & \frac{8 \alpha_1 }{9} & \frac{4 \alpha_1 }{3} & \frac{8 \alpha_1 }{9} \\
 -\frac{8 \alpha_1 ^2}{81}-\frac{\alpha_1  u_m}{3} & \frac{\alpha_1 }{3} & \frac{199 \alpha_1 }{270}+u_m & -\frac{\alpha_1 }{6} & -\frac{22 \alpha_1 }{135} \\
 -\frac{5 \alpha_1  u_m}{9} & \frac{5 \alpha_1 }{9} & \frac{2 \alpha_1 }{5} & u_m & -\frac{2 \alpha_1 }{5} \\
 \frac{8 \alpha_1 ^2}{81}-\frac{\alpha_1  u_m}{3} & \frac{\alpha_1 }{3} & \frac{22 \alpha_1 }{135} & \frac{\alpha_1 }{6} & u_m-\frac{199 \alpha_1 }{270} \\
\end{pmatrix}
\end{equation}
with the shifted characteristic polynomial 

\begin{equation}
 60 c^5-c^3\left(83 \overline{\alpha_1} ^2+60\right)+23c\left(\overline{\alpha_1} ^4 +\overline{\alpha_1} ^2\right) .
\end{equation}
The real eigenvalues of the HSSWME-L3 take the simple form
    \begin{equation}
            \lambda_{1,\ldots,5} \in \left\{ 
            u_m \pm \sqrt{\frac{23}{60}} \alpha_1 ,u_m,u_m \pm \sqrt{\alpha_1 ^2+gh}
 \right\}.
\end{equation}

\subsubsection{(H)SSWME-Q2}
The SSWME-Q2 system matrix
    \eqref{eq:ASysSSWME-Q2} yields the shifted characteristic polynomial
    \begin{multline}
        8960 c^4-9600 c^3 \left(\overline{s_1}-\overline{s_2}\right)-32 c^2 \left(783 \overline{s_1}{}^2+1458 \overline{s_2} \overline{s_1}+783 \overline{s_2}{}^2+280\right)\\
        -24 c \left(9 \overline{s_1}{}^3+405 \overline{s_2} \overline{s_1}{}^2-405 \overline{s_2}{}^2 \overline{s_1}-400 \overline{s_1}-9 \overline{s_2}{}^3+400 \overline{s_2}\right)\\
        +9 \left(783 \overline{s_1}{}^4+4068 \overline{s_2} \overline{s_1}{}^3+6426 \overline{s_2}{}^2 \overline{s_1}{}^2+208 \overline{s_1}{}^2+4068 \overline{s_2}{}^3 \overline{s_1}+1376 \overline{s_2} \overline{s_1}+783 \overline{s_2}{}^4+208 \overline{s_2}{}^2\right)=0.\\
    \end{multline}
    
    The restriction $s_1=s_2=: \frac{\alpha_1}{3}$ leads to the hyperbolic HSSWME-Q2 system matrix \eqref{eq:ASysHSSWME-Q2} 
with real eigenvalues 
    \begin{equation}
    \lambda_{1,\ldots,4} \in \{u_m \pm \frac{2 \alpha_1}{\sqrt{5}},u_m\pm\sqrt{gh+4 \alpha_1^2} \}.
\end{equation}  

\subsubsection{(H)SSWME-Q3}
    The SSWME-Q3 system matrix is given by $A_{sys}^{Q3} =$
    \begin{multline}
        \left( \begin{array}{c c c :}
 0 & 1 & 0 \\
 gh-u_m^2-\frac{48 s_1^2}{25}-\frac{156 s_2 s_1}{125}+\frac{96 s_3 s_1}{125}-\frac{204 s_2^2}{125}-\frac{48 s_3^2}{25}-\frac{156 s_2 s_3}{125} & 2 u_m & \frac{96 s_1}{25}+\frac{156 s_2}{125}-\frac{96 s_3}{125} \\
 -2 s_1 u_m-\frac{2161 s_1^2}{875}-\frac{8987 s_2 s_1}{3500}-\frac{2 s_3 s_1}{5}-\frac{17 s_2^2}{70}+\frac{611 s_3^2}{875}+\frac{2437 s_2 s_3}{3500} & 2 s_1 & u_m+\frac{5209 s_1}{1750}+\frac{4783 s_2}{7000}+\frac{127 s_3}{250} \\
 -2 s_2 u_m-\frac{96 s_1^2}{875}+\frac{492 s_2 s_1}{875}+\frac{96 s_3^2}{875}-\frac{492 s_2 s_3}{875} & 2 s_2 & \frac{1712 s_1}{875}+\frac{86 s_2}{875}-\frac{64 s_3}{125} \\
 -2 s_3 u_m-\frac{611 s_1^2}{875}-\frac{2437 s_2 s_1}{3500}+\frac{2 s_3 s_1}{5}+\frac{17 s_2^2}{70}+\frac{2161 s_3^2}{875}+\frac{8987 s_2 s_3}{3500} & 2 s_3 & \frac{559 s_1}{1750}+\frac{2833 s_2}{7000}-\frac{23 s_3}{250} \\
\end{array} \right.\\ \left. 
\begin{array}{: c c}
 0 & 0 \\
 \frac{156 s_1}{125}+\frac{408 s_2}{125}+\frac{156 s_3}{125} & -\frac{96 s_1}{125}+\frac{156 s_2}{125}+\frac{96 s_3}{25} \\
 \frac{787 s_1}{700}-\frac{143 s_2}{280}+\frac{283 s_3}{700} & \frac{23 s_1}{250}-\frac{2833 s_2}{7000}-\frac{559 s_3}{1750} \\
 u_m+\frac{258 s_1}{175}-\frac{258 s_3}{175} & \frac{64 s_1}{125}-\frac{86 s_2}{875}-\frac{1712 s_3}{875} \\
 -\frac{283 s_1}{700}+\frac{143 s_2}{280}-\frac{787 s_3}{700} & u_m-\frac{127 s_1}{250}-\frac{4783 s_2}{7000}-\frac{5209 s_3}{1750} \\
\end{array} \right).
    \end{multline}
    The characteristic polynomial is left out for brevity. However, the regularized matrix for the linear profile with $s_1=\frac{\alpha_1}{8}, s_2=\frac{\alpha_1}{3}, s_3=\frac{\alpha_1}{8} $ is the new HSSWME-Q3 system matrix given by 
    \begin{equation}
        A_{sys,H}^{Q3} = \begin{pmatrix}
 0 & 1 & 0 & 0 & 0 \\
 -\frac{\alpha_1 ^2}{3}+gh-u_m^2 & 2 u_m & \frac{4 \alpha_1 }{5} & \frac{7 \alpha_1 }{5} & \frac{4 \alpha_1 }{5} \\
 -\frac{5 \alpha_1 ^2}{36}-\frac{\alpha_1  u_m}{4} & \frac{\alpha_1 }{4} & \frac{199 \alpha_1 }{300}+u_m & \frac{\alpha_1 }{48} & -\frac{49 \alpha_1 }{300} \\
 -\frac{2 \alpha_1  u_m}{3} & \frac{2 \alpha_1 }{3} & \frac{16 \alpha_1 }{75} & u_m & -\frac{16 \alpha_1 }{75} \\
 \frac{5 \alpha_1 ^2}{36}-\frac{\alpha_1  u_m}{4} & \frac{\alpha_1 }{4} & \frac{49 \alpha_1 }{300} & -\frac{\alpha_1 }{48} & u_m-\frac{199 \alpha_1 }{300} \\
\end{pmatrix}
    \end{equation}
with real eigenvalues 
    \begin{equation}
    \lambda_{1,\ldots,5} \in \{u_m\pm\sqrt{\frac{19}{45}\alpha_1^2},u_m,u_m\pm\sqrt{gh+\alpha_1^2}\}.
\end{equation} 
\end{section}
\printbibliography
\end{document}